\documentclass[11pt]{article}
\usepackage{amsfonts}
\usepackage{eucal,color}
\usepackage{amssymb,epsfig,latexsym}
\usepackage{graphicx,amsmath,amssymb,latexsym,psfrag}
\usepackage{graphics,graphicx}
\newtheorem{Theorem}{Theorem}[section]
\newtheorem{Proposition}[Theorem]{Proposition}

\newtheorem{Definition}[Theorem]{Definition}

\newtheorem{Lemma}[Theorem]{Lemma}
\newtheorem{Remark}[Theorem]{Remark}

\newcommand{\re}{\mathrm{e}}
\newcommand{\im}{\mathbf i}

\begin{document}
\title{
Asymptotic analysis of Poisson shot noise processes, and applications}
\author{Giovanni Luca Torrisi\thanks{Consiglio Nazionale delle Ricerche, Via dei Taurini 19, 00185 Roma, Italy.
e-mail: \tt{giovanniluca.torrisi@cnr.it}}\,\,\,and Emilio Leonardi\thanks{Dipartimento di Elettronica e Telecomunicazioni,
Politecnico di Torino, Italy. e-mail: \tt{emilio.leonardi@polito.it}}}
\date{}
\maketitle
\begin{abstract}
Poisson shot noise processes are natural generalizations of compound Poisson processes that have been widely applied in insurance, neuroscience, seismology, computer science and epidemiology. In this paper we study sharp deviations, fluctuations and the stable probability approximation of Poisson shot noise processes. Our achievements extend, improve and complement existing results in the literature. 
We apply the theoretical results to Poisson cluster point processes, including generalized linear Hawkes processes, and risk processes with delayed claims.
Many examples are discussed in detail.
\end{abstract}

\section{Introduction}

We consider Poisson shot noise processes $\{S_t\}_{t>0}$ of the form
\begin{equation}\label{eq:PSN}
S_t:=\sum_{n\geq 1}H(t-T_n,M_n)\mathbf{1}_{(0,t]}(T_n),
\end{equation}
where $\{T_n\}_{n\geq 1}$ is a homogeneous Poisson process on $(0,\infty)$ with intensity $\lambda>0$, $\{M_n\}_{n\geq 1}$ is a 
sequence of random variables with values on some measurable space $(\mathrm{M},\mathcal{M})$, independent of $\{T_n\}_{n\geq 1}$, $H:[0,\infty)\times\mathrm{M}\to\mathbb R$ is a measurable function and $\bold{1}_A(\cdot)$ is the indicator function of a set $A$. We suppose that the random variables $\{M_n\}_{n\geq 1}$ are independent and identically distributed. 

Poisson shot noise processes are natural generalizations of compound Poisson processes, which
have found applications in different fields, due to their versatility and mathematical tractability. In insurance mathematics, Poisson shot noise processes arise as models of Incurred But Not Reported Claims \cite{KM0, KM}. In this context, $T_n$ represents the instant at which the $n$th claim arrives, 
for any $m\in\mathrm M$, the function $H(\cdot,m)$ is non-negative and non-decreasing, the random quantity $H(\infty,M_n)$ models the total pay-off caused by the $n$th claim and the random function $H(\cdot-T_n,M_n)$ models the evolution of the pay-off process for the $n$th claim. We refer the reader to \cite{B, GT1, GMT2, KM0, KM, MT, MST, T} (and the literature cited therein) for specific applications in insurance mathematics of Poisson shot noise processes, such as estimates of ruin probabilities of risk processes with delay in claim settlement.  In neurophysiology, Poisson shot noise processes appear as models of synaptic input. In this context,  $T_n$ describes the $n$th presynaptic event happened in the time interval $(0,t]$, $H(\cdot,\cdot)$ is the impulse-response function and $M_n$ models a possible synaptic inhomogeneity \cite{BD, P}. In computer science, Poisson shot noise processes are used e.g. as traffic, queues or caches models \cite{BTIDO,GMT1, KL, K, LT, T}.  Poisson shot noise processes are also exploited to model earthquake aftershocks \cite{VJ} and epidemics \cite{M}.
Spatial versions of $\{S_t\}_{t>0}$ are proposed in \cite{BB} to model the interference in wireless communication networks, see also \cite{BBbook,GT2, PT, TL}.  We study sharp estimates, fluctuations and stable probability
approximation of spatial Poisson shot noise processes, with applications to communication networks, in a companion paper.

Over the years, the mathematics of Poisson shot noise processes have been investigated by many authors. The central limit theorem and the Berry-Esseen bound have been proved in \cite{LA1,LA2}. The law 
of the large numbers and functional central limit theorems have been studied in \cite{KM}. Scalar and sample path large deviations are investigated in \cite{B, GMT1, MT, ST}.
In \cite{KMS}, the authors study the weak convergence to a multivariate infinite-variance stable distribution of the finite-dimensional distributions of a properly normalized and centred Poisson shot noise process.

In words, the main theoretical contributions of this paper are: $(i)$ The sharp deviations, at scales $O(t)$, of $\{S_t\}_{t>0}$ from its asymptotic mean; this result improves the tail estimates based on large deviations given in \cite{B,MT}; $(ii)$ The fluctuations, at scales $o(t)$, of $\{S_t\}_{t>0}$ from its asymptotic mean; these results improve the central limit theorem in \cite{LA1}; $(iii)$ In the case of a multiplicative noise,  we provide quantitative limit theorems for the weak convergence, as $t\to\infty$, of
$S_t$ (properly re-scaled) to a random variable $S$ with a stable law; these results complement the research started in \cite{KMS}. 
We emphasize that the results about stable approximations of Poisson shot noise
processes cover only stable laws with stability parameter $\alpha$ and skewness parameter $\beta$ such that either $\alpha\neq 1$ or $\alpha=1$ and $\beta=0$
(see Subsection \ref{subsec:stable}).  We remark that, although stable laws
have been defined also for $\alpha=1$ and $\beta\neq 0$, the most relevant and well-known  stable laws (i.e., the Gaussian, the Cauchy and the L\'evy) are all encompassed by our study. 

From the point of view of the applications, our main achievements concern: $(i)$ The sharp deviations and fluctuations of Poisson cluster processes
and, in particular, of generalized linear Hawkes processes,  extending in this way the results in \cite{GZ}; $(ii)$ Estimates of ruin probabilities
of risk processes with delayed claims,  refining the large deviation approximation provided in \cite{B}. 

Most of our theoretical contributions are achieved by means of the recently developed mod-$\phi$ convergence theory \cite{FMN, FMN2}.  Roughly speaking, 
provided that  a natural normalization of the characteristic function of a stochastic process converges to some non-trivial limiting function,   mod-$\phi$ convergence theory allows us 
to obtain  precise deviations  and fluctuations of the process from its asymptotic mean, improving classical results stemming from  large deviations and central limit theorems.

The paper is structured as follows.  In Section \ref{sec:modfi} we give some preliminaries on mod-$\phi$ convergence theory,  compound Poisson and stable laws. Moreover, we state an elementary inequality between complex numbers and the Fa\`a di Bruno formula which will be exploited a lot of times throughout the paper.
In Section \ref{sec:mainteo} we present the results on sharp deviations and fluctuations from the asymptotic mean of Poisson shot noise processes.  Applications to Poisson cluster processes and ruin probabilities are described in Sections \ref{sec:poissoncluster} and \ref{sec:insurance}, respectively. 
In particular,  in Section \ref{sec:poissoncluster} we shall consider extensions of the classical linear Hawkes process (see \cite{H}), where
the number of offspring of any parent is not necessarily Poisson distributed and the law of the birth times is not necessarily absolutely continuous with respect to the Lebesgue measure.  In Section \ref{sec:stablegrosso} we state the results concerning the stable approximation of Poisson shot noise processes (with a multiplicative noise), which extend
well-known results for the compound Poisson process (which are indeed recovered considering a constant shot shape).  All the proofs are given in Section \ref{sec:proofs}.

\section{Preliminaries}\label{sec:modfi}

\subsection{Mod-$\phi$ convergence}

We preliminary recall that a real-valued random variable $X$ (or its law) is said infinitely divisible if, for any $n\in\mathbb N:=\{1,2,\ldots\}$,
\[
X\overset{d}{=}X_1+\ldots+X_n,
\]
for some independent and identically distributed random variables $X_1,\ldots,X_n$. Here the symbol $\overset{d}{=}$ denotes the equality in law.

We proceed providing the definition of mod-$\phi$ convergence, 
see \cite{FMN}.

\begin{Definition}\label{def:mod}
Let $\phi$ be a non-constant infinitely divisible law on $\mathbb R$ and let $D\subseteq\mathbb C$ be a subset of the complex plane which contains $0$. We assume that the Laplace transform of $\phi$ is defined on $D$, i.e.,
\[
\Big|\int_{\mathbb R}\mathrm{e}^{z x}\,\phi(\mathrm{d}x)\Big|<\infty,\quad\text{for all $z\in D$,}
\]
and it has L\'evy exponent $\eta(\cdot)$ on $D$, i.e.,
\[
\int_{\mathbb R}\mathrm{e}^{z x}\,\phi(\mathrm{d}x)=\mathrm{e}^{\eta(z)},\quad\text{$z\in D$.}
\]
Let $\{X_t\}_{t>0}$ be a real-valued stochastic process with Laplace transform defined on $D$, i.e.,
such that $|\mathbb{E}[\mathrm{e}^{z X_t}]|<\infty$, for all $t>0$ and $z\in D$, and let $x(t)$, $t>0$, be a positive function such that
$x(t)\to+\infty$, as $t\to+\infty$.

We say that $\{X_t\}_{t>0}$ converges mod-$\phi$ on $D$,
with parameter function $x(\cdot)$ and limiting function $\psi(\cdot)$, if $\psi:D\to\mathbb C$ is analytic, it does not vanish on $\mathrm{Re}D:=\{\mathrm{Re}z:\,\,z\in D\}$ and 
\begin{equation*}
\text{for any compact $K\subseteq D$,}\quad\lim_{t\to\infty}\sup_{z\in K}|\psi_t(z)-\psi(z)|=0,
\end{equation*}
where 
\begin{equation*}
\psi_t(z):=\mathbb{E}[\mathrm{e}^{z X_t}]\mathrm{e}^{-x(t)\eta(z)}.
\end{equation*}

We say that $\{X_t\}_{t>0}$ converges mod-$\phi$ on $D$ with speed $O(x(t)^{-\sigma})$, for some $\sigma\in\mathbb N$,
and limiting function $\psi(\cdot)$, if $\psi:D\to\mathbb C$ is analytic, it does not vanish on $\mathrm{Re}D$ and 
\begin{equation*}
\text{for any compact $K\subseteq D$ there exists $C_K>0$:}\quad\sup_{z\in K}|\psi_t(z)-\psi(z)|\leq C_K x(t)^{-\sigma}.
\end{equation*}
\end{Definition}

In this paper we consider two different classes of infinitely divisible (reference) laws:
the compound Poisson law and the stable law.

\subsection{Compound Poisson laws}

Let $X$ be a non-negative random variable such that
\begin{equation}\label{eq:compLT}
a_X\in (0,\infty],\quad\text{where $a_X:=\sup\{\gamma:\,\,\mathbb{E}[\mathrm{e}^{\gamma X}]<\infty$\}.}
\end{equation}
Throughout this paper,  we denote by
$\phi_{\lambda,X}$ the compound Poisson law with L\'evy exponent
\begin{equation}\label{eq:levyexpintro9Ago}
\eta_{\lambda,X}(z):=\lambda(\mathbb{E}[\mathrm{e}^{zX}]-1),\quad z\in D_{cp}(a_X):=\{z\in\mathbb C:\,\,\mathrm{Re}z<a_X\},
\end{equation}
i.e.,  a random variable $Y$ has law $\phi_{\lambda,X}$ if and only if
\[
\mathbb{E}[\mathrm{e}^{z Y}]=\mathrm{e}^{\eta_{\lambda,X}(z)}, \quad z\in D_{cp}(a_X).
\] 

\subsection{Stable laws}\label{subsec:stable}

Let $c>0$ (the scale parameter), $\alpha\in (0,2]$ (the stability parameter) and $\beta\in [-1,1]$
(the skewness parameter) be fixed.
We consider the stable law $\phi_{c,\alpha,\beta}$ with parameters $(c,\alpha,\beta)$, whose Fourier transform
(or characteristic function) has L\'evy exponent
\begin{equation*}
\eta_{c,\alpha,\beta}(\bold i\xi):=-|c\xi|^{\alpha}(1-\bold{i}\beta h(\alpha,\xi)\mathrm{sgn}(\xi)),\quad\xi\in\mathbb{R}\setminus\{0\},
\quad\eta_{c,\alpha,\beta}(0):=0,
\end{equation*}
i.e., a random variable $S$ has law $\phi_{c,\alpha,\beta}$ if and only if $\mathbb{E}[\mathrm{e}^{\bold{i}\xi S}]=\mathrm{e}^{\eta_{c,\alpha,\beta}(\bold i\xi)}$, $\xi\in\mathbb R$.
Here
\[
h(\alpha,\xi):=\bold{1}\{\alpha\neq 1\}\tan\left(\frac{\pi\alpha}{2}\right)-\bold{1}\{\alpha=1\}\frac{2}{\pi}\log|\xi|\quad\text{and}\quad
\text{$\mathrm{sgn}(\xi)$ denotes the sign of $\xi$.}
\]
Since 
\begin{equation}\label{eq:moduluscarfunc}
|\mathrm{e}^{\eta_{c,\alpha,\beta}(\bold i \xi)}|=\mathrm{e}^{-|c\xi|^{\alpha}},\quad\xi\in\mathbb R,
\end{equation}
the characteristic function of $\phi_{c,\alpha,\beta}$ is integrable. Therefore, the stable law 
with parameters $(c,\alpha,\beta)$ has a density with respect to the Lebesgue measure. A standard computation shows that such a density
is bounded above by $\frac{\Gamma\left(\frac{1}{\alpha}\right)}{\alpha\pi c}$, where $\Gamma(\cdot)$ denotes the Euler gamma function. 

For later purposes, we recall the following scaling property
of the L\'evy exponent $\eta_{c,\alpha,\beta}(\cdot)$. For any $t>0$ and $\xi\in\mathbb R$, it holds:
\begin{equation}\label{eq:scale1}
t\eta_{c,\alpha,\beta}\left(\frac{\bold i\xi}{t^{1/\alpha}}\right)=\eta_{c,\alpha,\beta}(\bold i\xi)\quad\text{if either $\alpha\neq 1$ or $\alpha=1$ and $\beta=0$.}
\end{equation}
Finally, we recall some famous stable laws: the standard Gaussian distribution corresponds to the stable law with parameters $(2^{-1/2},2,0)$, the standard Cauchy distribution corresponds to the stable law
with parameters $(1,1,0)$ and the standard L\'evy distribution corresponds to the stable law with parameters $(1,2^{-1},1)$.

\subsection{An elementary inequality and the Fa\`a di Bruno formula}

Throughout this paper we exploit the following elementary inequality between complex numbers:

\begin{Lemma}\label{le:standardineq}
It holds:
\begin{equation*}
|\mathrm{e}^{z_1}-\mathrm{e}^{z_2}|\leq |z_1-z_2|\mathrm{e}^{\max\{\mathrm{Re}z_1,\mathrm{Re}z_2\}},\quad z_1,z_2\in\mathbb{C}.
\end{equation*}
\end{Lemma}
Since we have not found a proof of this inequality in standard textbooks of complex analysis, we show it in Section \ref{sec:proofs}.

Hereafter, for a sufficiently smooth function $f$ we denote by $f^{(n)}$
its derivative of order $n\in\mathbb N$.

\begin{Lemma}\label{le:FaadiBruno} $($Fa\`a\,\,di\,\,Bruno\,\,formula$)$
For any sufficiently smooth functions $g$ and $h$,
\[
(g\circ h)^{(j)}(x)=j!\sum_{i=1}^{j}\frac{g^{(i)}(h(x))}{i!}\sum_{m_1+m_2+\ldots+m_i=j}\frac{h^{(m_1)}(x)}{m_1!}\ldots\frac{h^{(m_i)}(x)}{m_i!},
\quad j\in\mathbb N,
\]
where the sum is taken over all the $m_1,\ldots,m_i\in\mathbb N$ such that $m_1+\ldots+m_i=j$. 
\end{Lemma}

\section{Sharp deviations and fluctuations of Poisson shot noise processes}\label{sec:mainteo}

As already mentioned in the Introduction,  our analysis relies on the mod-$\phi$ convergence theory.  Specifically,  we obtain different sharp deviations estimates  (at scale $O(t)$) of Poisson shot noise processes depending on whether the reference measure $\phi$ is non-lattice or lattice,  as it can be realized by 
comparing the formulas \eqref{eq:deviationordert} and \eqref{eq:latticemaggiugbreve}. 
The distinction \lq\lq $\phi$ non-lattice" and \lq\lq $\phi$ lattice", instead, has no impact on the results about the fluctuations (at scales $o(t)$) of Poisson shot noise processes.
For the sake of completeness, we recall that a probability law is lattice if its support is included in a set of the form $\gamma_1+\gamma_2\mathbb{Z}$, for some parameters $\gamma_1\in\mathbb{R}$ and $\gamma_2>0$. 

In both cases (non-lattice and lattice), the shot shape $H(\cdot,\cdot)$ is supposed to be a non-negative function.
Moreover we   assume  that
\begin{equation}\label{eq:lighttailCP}
\text{$Z:=\sup_{t\geq 0}H(t,M_1)$ is such that $a:=a_Z\in (0,\infty]$.} 
\end{equation}
Although in the non-lattice case,  i.e.,  under the assumption
\begin{equation}\label{eq:nonlattice}
\text{$\phi_{\lambda,Z}$ is non-lattice,} 
\end{equation}
we simply suppose that:
\begin{equation}\label{eq:noneg}
\text{function $H:[0,\infty)\times\mathrm{M}\to [0,\infty)$ is non-negative,}
\end{equation}
in the lattice case,  we naturally assume that:
\begin{equation}\label{eq:integerpos}
\text{function $H:[0,\infty)\times\mathrm{M}\to\mathbb{N}\cup\{0\}$ takes values in $\mathbb{N}\cup\{0\}$.} 
\end{equation}
A crucial hypothesis to prove the mod-$\phi_{\lambda,Z}$ convergence of Poisson shot noise processes is
\begin{equation}\label{eq:Lq}
\text{$\int_0^\infty(Z-H(s,M_1))\,\mathrm{d}s\in L^q(\mathbb P)$, for any $q>1$.}
\end{equation}
However,  
when sharp deviations are concerned and $\phi_{\lambda,Z}$ is lattice we need to strengthen \eqref{eq:Lq} by assuming
\begin{equation}\label{eq:Lqspeed}
\text{ 
$\exists$ $\sigma\in\mathbb N$ and $\{\kappa_q\}_{q>1}\subset (0,\infty)$:}
\end{equation}
\begin{equation*}
\text{$\sup_{q>1}\kappa_q^{-1}\Big\|\int_t^\infty(Z-H(s,M_1))\,\mathrm{d}s \Big\|_{L^q(\mathbb P)}\leq t^{-\sigma}$, for all $t$ large enough.}
\end{equation*}
Indeed, this condition guarantees the mod-$\phi_{\lambda,Z}$ convergence of Poisson shot noise processes with speed $O(t^{-\sigma})$. Finally, we mention that to prove the fluctuations of Poisson shot noise processes when $\phi_{\lambda,Z}$ is non-lattice, we need to assume that $\phi_{\lambda,Z}$ has a density, i.e., 
\begin{equation}\label{eq:absolutecontinuous}
\text{$\phi_{\lambda,Z}$ is absolutely continuous with respect to the Lebesgue measure.} 
\end{equation}

Hereon,  under \eqref{eq:lighttailCP} (and either \eqref{eq:noneg} or \eqref{eq:integerpos})
\begin{equation*}
\text{$\forall$
$x\in (0,\lambda\mathbb{E}[Z\mathrm{e}^{a Z}])$, we denote by $\theta_x\in (-\infty,a)$ the unique solution of the equation $\lambda\mathbb{E}[Z\mathrm{e}^{\theta Z}]=x$.}
\end{equation*}
Note that $0<x<\lambda\mathbb{E}[Z]$ if and only if $\theta_x<0$; $x=\lambda\mathbb{E}[Z]$ if and only if $\theta_x=0$; 
$x\in(\lambda\mathbb{E}[Z],\lambda\mathbb{E}[Z\mathrm{e}^{a Z}])$ if and only if $\theta_x\in (0,a)$. 

Throughout this article, we denote by $\eta_{\lambda,Z}^*$
the Fenchel-Legendre transform of $\eta_{\lambda,Z}$ (see \eqref{eq:levyexpintro9Ago}), i.e.,
\[
\eta_{\lambda,Z}^*(x):=\sup_{\theta\in\mathbb R}(\theta x-\eta_{\lambda,Z}(\theta)),\quad x>0.
\]
Under \eqref{eq:lighttailCP} (and either \eqref{eq:noneg} or \eqref{eq:integerpos}), standard computations show that
\begin{equation}\label{eq:FenchelLegendre}
\eta_{\lambda,Z}^*(x)=x\theta_x-\lambda(\mathbb{E}[\mathrm{e}^{\theta_x Z}]-1),\quad x>0.
\end{equation}
Hereafter, for ease of notation, we also set
\begin{equation}\label{eq:varphi27Lug2021}
\varphi(z):=\int_0^\infty(\mathbb{E}[\mathrm{e}^{z H(s,M_1)}]-\mathbb{E}[\mathrm{e}^{z Z}])\,\mathrm{d}s,\quad z\in\mathbb C.
\end{equation}

\subsection{Sharp deviations at scales $O(t)$}

The following theorems provide the exact asymptotic behavior of the tail of the Poisson shot noise when the reference compound Poisson law $\phi_{\lambda,Z}$ is either
non-lattice or lattice, respectively.

\begin{Theorem}\label{thm:deviationsnonlattice} $($$\bold{Non}$-$\bold{lattice\,\,case}$$)$
Assume 
 \eqref{eq:lighttailCP}, \eqref{eq:nonlattice},  \eqref{eq:noneg} and \eqref{eq:Lq}. Then, for any $x\in (\lambda\mathbb{E}[Z],\lambda\mathbb{E}[Z\mathrm{e}^{a Z}])$, we have
\begin{equation}\label{eq:deviationordert}
\mathbb{P}(S_t\geq tx)=
\frac{\exp(-t\eta_{\lambda,Z}^*(x)+\lambda\varphi(\theta_x))}{\theta_x\sqrt{2\lambda\pi t\mathbb{E}[Z^2\mathrm{e}^{\theta_x Z}]}}
(1+o(1)),\quad\text{as $t\to+\infty$.}
\end{equation}

\end{Theorem}

\begin{Theorem}\label{thm:deviationslattice}$($$\bold{Lattice\,\,case}$$)$.
Assume
\eqref{eq:lighttailCP},  \eqref{eq:integerpos}  and \eqref{eq:Lqspeed}. 
Then:\\
\noindent$(i)$ For any $x\in(0,\lambda\mathbb{E}[Z\mathrm{e}^{a Z}])$ such that $tx\in\mathbb N$, we have
\begin{equation}\label{eq:latticeugbreve}
\mathbb{P}(S_t=tx)=
\frac{\exp(-t\eta_{\lambda,Z}^*(x)+\lambda\varphi(\theta_x))}{\sqrt{2\lambda\pi t\mathbb{E}[Z^2\mathrm{e}^{\theta_x Z}]}}
(1+o(1)),\quad\text{as $t\to+\infty$.}
\end{equation}
More in general,  for any $x\in(0,\lambda\mathbb{E}[Z\mathrm{e}^{a Z}])$ such that $tx\in\mathbb N$,  we have
\begin{align}
\mathbb{P}(S_t=tx)=
\frac{\exp(-t\eta_{\lambda,Z}^*(x))}{\sqrt{2\lambda\pi t\mathbb{E}[Z^2\mathrm{e}^{\theta_x Z}]}}&\Biggl(
\mathrm{e}^{\lambda\varphi(\theta_x)}
+\sum_{k=1}^{\sigma-1}\frac{a_{k}(\theta_x)}{t^k}+O\left(\frac{1}{t^\sigma}\right)\Biggr),\quad\text{as $t\to+\infty$,}\label{eq:deviationordertlattice}
\end{align}
where the quantities $a_k(\theta_x)$ are computed in Proposition \ref{prop:akbk}. \\
\noindent$(ii)$ For any $x\in(\lambda\mathbb{E}[Z],\lambda\mathbb{E}[Z\mathrm{e}^{a Z}])$ such that $tx\in\mathbb{N}$,  we have
\begin{equation}\label{eq:latticemaggiugbreve}
\mathbb{P}(S_t\geq tx)=
\frac{\exp(-t\eta_{\lambda,Z}^*(x)+\lambda\varphi(\theta_x))}{\sqrt{2\lambda\pi t\mathbb{E}[Z^2\mathrm{e}^{\theta_x Z}]}}\frac{1}{1-\mathrm{e}^{-\theta_x}}
(1+o(1)),\quad\text{as $t\to+\infty$.}
\end{equation}
More in general,  for any $x\in(\lambda\mathbb{E}[Z],\lambda\mathbb{E}[Z\mathrm{e}^{a Z}])$ such that $tx\in\mathbb{N}$,  we have
\begin{align}
\mathbb{P}(S_t\geq tx)=
\frac{\exp(-t\eta_{\lambda,Z}^*(x))}{\sqrt{2\lambda\pi t\mathbb{E}[Z^2\mathrm{e}^{\theta_x Z}]}}&
\frac{1}{1-\mathrm{e}^{-\theta_x}}
\Biggl(\mathrm{e}^{\lambda\varphi(\theta_x)}
+\sum_{k=1}^{\sigma-1}\frac{b_{k}(\theta_x)}{t^k}+O\left(\frac{1}{t^\sigma}\right)\Biggr),\,\,\text{as $t\to+\infty$,}\label{eq:deviationordertlatticetail}
\end{align}
where the quantities $b_k(\theta_x)$ are computed in Proposition \ref{prop:akbk}. 
\end{Theorem}

We emphasize that Theorem \ref{thm:deviationsnonlattice} and Theorem \ref{thm:deviationslattice}$(ii)$ refine the tail estimates provided by the corresponding large deviations
principle in \cite{GMT1} (see also \cite{MT}).

We shall give many examples where the quantities $\theta_x$, $\mathbb{E}[\mathrm{e}^{\theta_x Z}]$ and 
$\mathbb{E}[Z^2\mathrm{e}^{\theta_x Z}]$ can be computed. We shall give also some examples where the integral $\varphi(\theta_x)$
can be computed.  When this is not possible, our starting point to provide estimates of the integral is the following elementary proposition.

\begin{Proposition}\label{prop:estimateintegral}
Under the foregoing assumptions and notation, we have:\\
\noindent$(i)$ If $\theta_x<0$ $($i.e., $x\in (0,\lambda\mathbb{E}[Z])$$)$, then
\[
0\leq\varphi(\theta_x)
\leq-\theta_x\int_0^\infty\mathbb{E}[Z-H(s,M_1)]\,\mathrm{d}s.
\]
\noindent$(ii)$ If $\theta_x>0$ $($i.e., $x\in (\lambda\mathbb{E}[Z],\mathbb{E}[Z\mathrm{e}^{a Z}])$$)$, then for any $q,q'>1$ such that
$q^{-1}+q'^{-1}=1$, we have
\[
-\theta_x\|\mathrm{e}^{\theta_x Z}\|_{L^q(\mathbb P)}
\Big\|\int_0^\infty(Z-H(s,M_1))\,\mathrm{d}s\Big\|_{L^{q'}(\mathbb P)}
\leq\varphi(\theta_x)\leq 0.
\]
\end{Proposition}

We remark that this proposition is exploited e.g.  in the application to Poisson cluster processes to exhibit explicit estimates of 
the integral $\varphi(\theta_x)$, see Proposition \ref{prop:PoissonCluster}$(iii)$ and Proposition \ref{prop:GW}$(vi)$.

\subsection{Fluctuations at scales $o(t)$}

The fluctuations of the Poisson shot noise process from its asymptotic mean, at any scale which is $o(t)$ as $t\to\infty$, are provided by the following theorem, where $\mathrm{N}(0,1)$
denotes a random variable distributed according to the standard Gaussian law. 
We refer the reader to Remark \ref{re:scalings} for a brief discussion on the range of the scaling functions for the process $\{S_t-\lambda\mathbb{E}[Z]t\}_{t>0}$.

\begin{Theorem}\label{thm:fluctuations}
Assume either
 \eqref{eq:lighttailCP}, \eqref{eq:noneg}, \eqref{eq:Lq} and 
\eqref{eq:absolutecontinuous} or  \eqref{eq:lighttailCP},  \eqref{eq:integerpos} and \eqref{eq:Lq}. Then:\\
\noindent$(i)$ For any function $y(\cdot)$ such that $y(t)=o(t^{1/6})$, it holds
\begin{equation}\label{eq:CLT}
\mathbb{P}\left(\frac{S_t-\lambda\mathbb{E}[Z]t}{\sqrt{\lambda\mathbb{E}[Z^2]t}}>y(t)\right)=\mathbb{P}(\mathrm{N}(0,1)> y(t))
(1+o(1)),\quad\text{as $t\to+\infty$.}
\end{equation}
In particular, if $y(t)\to\infty$ as $t\to+\infty$, then $($by the asymptotics of the tail of $\mathrm{N}(0,1)$$)$
\begin{equation*}
\mathbb{P}\left(\frac{S_t-\lambda\mathbb{E}[Z]t}{\sqrt{\lambda\mathbb{E}[Z^2]t}}>y(t)\right)=\frac{\mathrm{e}^{-\frac{y(t)^2}{2}}}{y(t)\sqrt{2\pi}}
(1+o(1)),\quad\text{as $t\to+\infty$.}
\end{equation*}
\noindent$(ii)$ For any function $y(\cdot)$ such that $y(t)\to+\infty$, as $t\to\infty$, and $y(t)=o(t^{1/2})$, it holds
\begin{align}
\mathbb{P}\left(\frac{S_t-\lambda\mathbb{E}[Z]t}{\sqrt{\lambda\mathbb{E}[Z^2]t}}>y(t)\right)&=
\frac{\mathrm{e}^{-t\eta_{\lambda,Z}^*(v(t))}}
{\theta_{v(t)}\sqrt{2\lambda\pi t\mathbb{E}[Z^2\mathrm{e}^{\theta_{v(t)}Z}]}}(1+o(1))\nonumber\\
&=\frac{\mathrm{e}^{-t\eta_{\lambda,Z}^*(v(t))}}{y(t)\sqrt{2\pi}}(1+o(1)),\quad\text{as $t\to+\infty$,}\label{eq:extCLT}
\end{align}
and
\begin{equation}\label{eq:Dic3mattino1}
t\,\eta_{\lambda,Z}^*(v(t))=\frac{y(t)^2}{2}(1+o(1)).
\end{equation}
Here
\[
v(t):=\lambda\mathbb{E}[Z]+\frac{y(t)}{\sqrt t}\sqrt{\lambda\mathbb{E}[Z^2]}.
\]
\noindent$(iii)$ For any function $y(\cdot)$ such that $y(t)\to+\infty$, as $t\to\infty$, and $y(t)=o\left(t^{\frac{1}{2}-\frac{1}{m}}\right)$, $m\geq 3$ integer, it holds
\begin{align}
\mathbb{P}\left(\frac{S_t-\lambda\mathbb{E}[Z]t}{\sqrt{\lambda\mathbb{E}[Z^2]t}}>y(t)\right)
&=\frac{\exp\left(-\frac{y(t)^2}{2}\left(1+2\sum_{j=1}^{m-2}(\lambda\mathbb{E}[Z^2])^{(j+2)/2}
\frac{\theta_{\lambda\mathbb{E}[Z]}^{(j+1)}}{(j+2)!}\left(\frac{y(t)}{\sqrt t}\right)^j\right)\right)
}{y(t)\sqrt{2\pi}}(1+o(1)),\label{eq:fluctuationderivatives}
\end{align}
as $t\to+\infty$.
Here, the derivatives of $\theta_x$ evaluated at $\lambda\mathbb{E}[Z]$, denoted by $\theta_{\lambda\mathbb{E}[Z]}^{(j+1)}$, can be recursively computed by the formula:
\begin{equation}\label{eq:recursive}
\theta_{\lambda\mathbb{E}[Z]}^{(j+1)}=
j!\sum_{i=1}^{j}\frac{g^{(i)}(0)}{i!}\sum_{m_1+m_2+\ldots+m_i=j}
\frac{\theta_{\lambda\mathbb{E}[Z]}^{(m_1)}}{m_1!}\ldots\frac{\theta_{\lambda\mathbb{E}[Z]}^{(m_i)}}{m_i!},
\quad\text{$1\leq j\leq m-2$,}
\end{equation}
where the sum is taken over all the $m_1,\ldots,m_i\in\mathbb N$ such that $m_1+\ldots+m_i=j$,
\[
\theta'_{\lambda\mathbb{E}[Z]}=\frac{1}{\lambda\mathbb{E}[Z^2]}\quad\text{and}\quad g(x):=\frac{1}{\lambda\mathbb{E}[Z^2\mathrm{e}^{x Z}]}.
\]
\end{Theorem}

\subsection{Some remarks}

We conclude this section with the following three remarks.

\begin{Remark}\label{re:scalings}$($$\bold{On\,\,the\,\,range\,\,of\,\,the\,\,scaling\,\,functions}$$)$
Theorems \ref{thm:deviationsnonlattice} and \ref{thm:fluctuations} cover the whole range of scalings for the process $\{S_t-\lambda\mathbb{E}[Z]t\}_{t>0}$ up to the order of $t$. Indeed, Theorem \ref{thm:fluctuations}$(i)$ covers scalings, say $s(t)$, of $S_t-\lambda\mathbb{E}[Z]t$ such that 
$s(t)\ll t^{1/6}\sqrt t=t^{2/3}$. 
Theorem \ref{thm:fluctuations}$(ii)$ covers scalings of 
$S_t-\lambda\mathbb{E}[Z]t$ such that either $s(t)\sim t^{1/6}\sqrt t=t^{2/3}$ or $t^{2/3}\ll s(t)\ll t$. Finally, 
Theorem \ref{thm:deviationsnonlattice} refers to scalings of $S_t-\lambda\mathbb{E}[Z]t$ of order $t$. 
\end{Remark}

\begin{Remark}\label{re:extendedCLT}$($$\bold{Central\,\,Limit\,\,Theorem\,\,and\,\,Extended\,\,Central\,\,Limit\,\,Theorem}$$)$
Let the assumptions of Theorem \ref{thm:fluctuations} prevail. Then it is well-known that
the classical Central Limit Theorem for the Poisson shot noise holds, i.e.,
\begin{equation*}
\frac{S_t-\mathbb{E}[S_t]}{\sqrt{\mathbb{V}\mathrm{ar}(S_t)}}\to\mathrm{N}(0,1)\quad\text{in law, as $t\to\infty$.}
\end{equation*}
$($see e.g. Theorem 2.3 in \cite{KM}$)$. 
This Central Limit Theorem 
can be retrieved by using Theorem \ref{thm:fluctuations}$(i)$. Indeed,
as we shall check later on $($see the comment after the statement of
Lemma \ref{thm:modphiR}$)$, the conditions \eqref{eq:lighttailCP}, \eqref{eq:noneg} and \eqref{eq:Lq}
guarantee 
\begin{equation}\label{eq:mediavarianza}
\mathbb{E}[S_t]=\lambda\mathbb{E}[Z]t+O(1)\quad\text{and}\quad\mathbb{V}ar(S_t)=\lambda\mathbb{E}[Z^2]t+O(1). 
\end{equation}
It follows
\begin{equation*}
\frac{S_t-\mathbb{E}[S_t]}{\sqrt{\mathbb{V}\mathrm{ar}(S_t)}}=\left(1+\frac{O(1)}{\sqrt t}\right)^{-1/2}\left(
\frac{S_t-\lambda\mathbb{E}[Z]t}{\sqrt{\lambda\mathbb{E}[Z^2]t}}
-\frac{O(1)}{\sqrt t}\right).
\end{equation*}
Therefore, for any $x\in\mathbb R$,
\begin{align}
\mathbb{P}\left(\frac{S_t-\mathbb{E}[S_t]}{\sqrt{\mathbb{V}\mathrm{ar}(S_t)}}>x\right)
&=\mathbb{P}\left(\frac{S_t-\lambda\mathbb{E}[Z]t}{\sqrt{\lambda\mathbb{E}[Z^2]t}}
>x(t)\right)
=\mathbb{P}(\mathrm{N}(0,1)>x(t))(1+o(1)),\label{eq:Dic3pom1}
\end{align}
where the latter relation is a consequence of \eqref{eq:CLT} since
\begin{equation}\label{eq:ytilde}
x(t):=x\left(1+\frac{O(1)}{\sqrt t}\right)^{1/2}+\frac{O(1)}{\sqrt t}=o(t^{1/6}).
\end{equation}
The Central Limit Theorem for the Poisson shot noise process stated at the beginning
follows combining
\eqref{eq:Dic3pom1} with the trivial relation
\[
\mathbb{P}(\mathrm{N}(0,1)>x(t))=\mathbb{P}(\mathrm{N}(0,1)>x)(1+o(1)).
\]
In fact, Theorem \ref{thm:fluctuations} is a big improvement of the classical Central Limit Theorem for the Poisson shot noise process, yielding the Extended Central Limit Theorem for the Poisson shot noise process.
\end{Remark}

\begin{Remark}\label{re:alternativecdz}
As it will be clear from the proofs, in all the theorems above, conditions \eqref{eq:Lq} and \eqref{eq:Lqspeed} can be replaced respectively by
\begin{equation*}
\text{$Z\in L^{\infty}(\mathbb P)$ and $\int_0^\infty(Z-H(s,M_1))\,\mathrm{d}s\in L^1(\mathbb P)$,}
\end{equation*}
and
\begin{equation*}
\text{$Z\in L^{\infty}(\mathbb P)$ and $\exists$ $\sigma\in\mathbb N$ and $\kappa>0$:}
\end{equation*}
\begin{equation*}
\text{$\Big\|\int_t^\infty(Z-H(s,M_1))\,\mathrm{d}s\Big\|_{L^1(\mathbb P)}\leq\kappa t^{-\sigma}$ for all $t$ large enough.}
\end{equation*}
Note also that if $Z\in L^\infty(\mathbb P)$, then $Z$ satisfies \eqref{eq:lighttailCP} with $a:=a_Z=+\infty$.
\end{Remark}

\section{Application to Poisson cluster processes}\label{sec:poissoncluster}

Let $N_n(\cdot)$, $n\geq 1$, be independent and identically distributed finite and simple point processes on $[0,\infty)$, 
independent on $\{T_n\}_{\geq 1}$.  We assume $N_n(\{0\})=1$ for any $n\geq 1$. Denoting by $\{S_{n,k}\}_{k\geq 0}$, $S_{n,0}:=0$ the (ordered) points of $N_n$, 
and interpreting $T_n$ as the ancestor of 
the point process $\theta_{T_n}N_n\equiv\{S_{n,k}+T_n\}_{k\geq 0}$ of offspring,  we have that, at time $t>0$, the total number of offspring 
generated by the ancestors (and including the ancestors) 
is equal to
\begin{align}
\sum_{n\geq 1}\bold{1}_{(0,t]}(T_n)\theta_{T_n}N_n([T_n,t])
&=\sum_{n\geq 1}\bold{1}_{(0,t]}(T_n)N_n([0,t-T_n])
=\sum_{n\geq 1}\bold{1}_{(0,t]}(T_n)H(t-T_n,M_n),\label{eq:poissoncluster}
\end{align}
where 
$M_n:=N_n$, i.e., $M_n$ is a random variable with values on $\mathrm{M}:=\mathcal{N}$,  i.e.,
the space of finite and simple counting measures on $[0,\infty)$
(endowed with the usual vague topology), and
$H:[0,\infty)\times\mathrm M\to\mathbb{N}\cup\{0\}$ is defined by
\[
H(t,m)=H(t,\mu):=\mu([0,t]).
\]
Thus, assuming that $\{T_n\}_{n\geq 1}$ is a homogeneous Poisson process with intensity $\lambda$,  at time $t>0$, the total number of offspring 
generated by the ancestors (including the ancestors) 
is equal to the Poisson shot noise process $S_t$ defined by \eqref{eq:poissoncluster}.

We put $L_n:=\sup_{k\geq 1}S_{n,k}$, $n\geq 1$,  and note that $L_n$ is the \lq\lq length" of the cluster point process $N_n(\cdot)$. 
Clearly, the random variables $\{L_n\}_{n\geq 1}$ are independent and identically distributed, and we set $L:=L_1$. 
The following proposition holds.

\begin{Proposition}\label{prop:PoissonCluster}
Let $\{S_t\}_{t>0}$ be the Poisson shot noise process \eqref{eq:poissoncluster}, i.e., $S_t$ denotes
the number of offspring generated by the Poissonian ancestors $\{T_n\}_{n\geq 1}$ in the time interval $(0,t]$ $($including the ancestors$)$.
Assume that $Z:=N_1([0,\infty))$ satisfies \eqref{eq:lighttailCP}, and that
\begin{equation}\label{eq:Lpiccolo}
\text{$\mathbb{E}[L^k]<\infty$, $\forall$ $k\in\mathbb N$.}
\end{equation}
Then, setting $H(s,M_1):=N_1([0,s])$, $s\in [0,t]$, we have that:\\
\noindent$(i)$ The formulas \eqref{eq:CLT}, \eqref{eq:extCLT} and
\eqref{eq:fluctuationderivatives} of Theorem \ref{thm:fluctuations} hold. \\
\noindent$(ii)$ For any arbitrarily fixed $\sigma\in\mathbb N$, the formulas \eqref{eq:deviationordertlattice} and \eqref{eq:deviationordertlatticetail} 
of Theorem \ref{thm:deviationslattice}  hold.\\
$(iii)$ The following estimates hold:\\
\noindent$(1)$ For any $x\in(0,\lambda\mathbb{E}[Z])$ and any $q,q'>1$ such that $q^{-1}+q'^{-1}=1$, we have
\begin{equation}\label{eq:lowerbd2800}
0\leq\varphi(\theta_x)
\leq-\theta_x\|L\|_{L^{q}(\mathbb P)}\|Z\|_{L^{q'}(\mathbb P)}.
\end{equation}
\noindent$(2)$ For any $x\in(\lambda\mathbb{E}[Z],\lambda\mathbb{E}[Z\mathrm{e}^{aZ}])$ and any $q,q',q_1,q_2>1$ such that
$q^{-1}+q'^{-1}=1$ and $q_1^{-1}+q_2^{-1}=1$, we have
\begin{equation}\label{eq:lowerbd2801}
-\theta_x\|\mathrm{e}^{\theta_x Z}\|_{L^q(\mathbb P)}\|L\|_{L^{q'q_1}(\mathbb P)}
\|Z\|_{L^{q'q_2}(\mathbb P)}
\leq\varphi(\theta_x)
\leq 0.
\end{equation}

Here the function $\varphi(\cdot)$ is given by \eqref{eq:varphi27Lug2021} with $Z$ and $H(\cdot,M_1)$ defined at the beginning of the statement.
\end{Proposition}

We note that the Part $(ii)$ of this proposition refines the tail estimates provided by the large deviations principle in \cite{BT}.

\subsection{Poisson cluster processes with a Galton-Watson branching structure or generalized linear Hawkes processes}\label{subsec:Kuno}

We consider Poisson cluster processes whose clusters have a Galton-Watson branching structure.
Clearly, to define such processes it suffices to describe the structure of the cluster point process $N_1$.

The points in the cluster $C_0$,  generated by a common ancestor placed in the origin, are partitioned in generations $g\in \mathbb{N}\cup\{0\}$.   Every point belonging to the $(g-1)$th generation will give rise to a random number of \lq\lq children\rq\rq\,  points in the $g$th generation.
By definition,  $N_1$ is the point process with support $C_{0}$. 
Let 
\begin{itemize}
	\item $K_{h}$ be  the number of points in the $h$th generation, 
	\item $W_{h}$ be  the total number of points in the cluster $C_{0}$ up the $h$th generation, 
	\item  $\{B_{1}^{h},\ldots,B_{K_{h}}^{h}\}$ be the birth times (i.e., the times at which points are placed) in the $h$th generation. 
\end{itemize}

We assume that  the birth times are arranged in increasing order, i.e.,
\[
B_{1}^{h}<\ldots<B_{K_h}^{h}.
\]

Let $\{P_{i,j}\}_{(i,j)\in\mathbb{N}^2}$ be a sequence 
of $\mathbb{N}\cup\{0\}$-valued independent random variables with law $\{p_k\}_{k\geq 0}$,  where $P_{ij}$ represents  the number of children generated by the $j$th individual of the $i$th generation. Let 
 $\{B_{i,j,k}\}_{(i,j,k)\in\mathbb{N}^3}$ be a sequence of non-negative independent and identically distributed random variables, 
 independent of $\{P_{i,j}\}_{(i,j)\in\mathbb{N}^2}$,  where $B_{ijk}$
 represents the time lag  between the birth time of the $j$th individual of the $i$th generation and the birth time 
of its $k$th child. 

The point $\{0\}$ constitutes the $0$th generation.   By construction we have  $K_0=1$ and  $W_0=1$. 
The birth time of the unique individual in this generation is given by $B_{1}^{0}:=0$.
The points in the  $g$th generation are generated by the points in the $(g-1)$th generation according to the following rule:
\\
\noindent $(i)$ If $K_{g-1}=0$, then the $(g-1)$th generation is empty and 
then the $g$th generation will be empty as well. We set $K_g:=0$ and $W_g:=W_{g-1}$. \\
\noindent $(ii)$ If $K_{g-1}>0$, then:
\begin{itemize}
\item $K_g:=P_{g,1}+\ldots+P_{g,K_{g-1}}$ is the number of births in the $g$th generation.
\item $W_g:=W_{g-1}+K_g$ is the total number of points in the cluster $C_{0}$ until the $g$th generation.
\item The birth times in the $g$th generation are given by the union of the following sets: 
\end{itemize}
\[
\{B_1^{g-1}+B_{g,1,1},\ldots,B_1^{g-1}+B_{g,1,P_{g,1}}\},
\]
which are the birth times of the children of the parent born at time $B_1^{g-1}$,
\[
\{B_2^{g-1}+B_{g,2,1},\ldots,B_2^{g-1}+B_{g,2,P_{g,2}}\},
\]
which are the birth times of the children of the parent born at time $B_2^{g-1}$
\[
\cdots
\]
\[
\{B_{K_{g-1}}^{g-1}+B_{g,K_{g-1},1},\ldots,B_{K_{g-1}}^{g-1}+B_{g,K_{g-1},P_{g,K_{g-1}}}\},
\]
which are the birth times of the children of the parent born at time $B_{K_{g-1}}^{g-1}$. 

The points in the $g$th 
generation are then arranged in increasing order as
\[
\{B_{1}^{g},\ldots,B_{K_{g}}^{g}\}. 
\]


If the the law of $B_{1,1,1}$ has a probability density $h(\cdot)/\int_0^\infty h(t)\,\mathrm{d}t$, where $h:(0,\infty)\to[0,\infty)$
is an integrable function, and the law of $P_{1,1}$ is Poisson with mean $\int_0^\infty h(t)\,\mathrm{d}t$, then
the corresponding Poisson shot noise process $\{S_t\}_{t>0}$, where $S_t$ denotes the number of points generated by the Poissonian ancestors
in the time interval $(0,t]$, is a classical Hawkes process (or linear Hawkes process) \cite{H}, for which large deviations were studied in \cite{BT}
and sharp deviations and fluctuations have been recently investigated in \cite{GZ}.

The sequence $\{K_n\}_{n\in\mathbb{N}\cup\{0\}}$ is a Galton-Watson process 
from an initial population of one individual and with offspring law $\{p_k\}_{k\geq 0}$. 
Hereon, we assume that 
\begin{equation}\label{eq:subcritico}
0<\mathbb{E}[P_{1,1}]<1,
\end{equation}
so that the Galton-Watson process is subcritical and the total progeny
\[
Z:=N_1([0,\infty))=\sum_{n\geq 0}K_n
\]
has mean equal to $(1-\mathbb{E}[P_{1,1}])^{-1}$. For later purposes, we recall that the distribution of the total progeny 
$Z:=N_1([0,\infty))$ is related to the offspring distribution by the formula
\begin{equation}\label{eq:totoff}
\mathbb{P}(Z=k)=\frac{1}{k}\mathbb{P}(P_{11}+\ldots+P_{1k}=k-1),\quad k\in\mathbb N
\end{equation}
(the reader is referred to e.g. \cite{Remco} for an introduction to branching processes). Hereafter, we also suppose that
\begin{equation}\label{hp:thetac}
b:=a_{P_{1,1}}\in (0,\infty],
\end{equation}
where $a_{P_{1,1}}$ is defined by \eqref{eq:compLT}. So (being the branching process subcritical) by Theorem 2.1 in \cite{SIG} we have that
the total progeny $Z$ satisfies \eqref{eq:lighttailCP}.

The \lq\lq length" of $N_1$ is given by
\[
L:=\sup_{g\in\mathbb N}\sup_{j\leq K_g}B_j^g.
\]
Define
\[
V_g:=\sum_{k=1}^{K_{g-1}}\sum_{j=1}^{P_{g,k}}B_{g,k,j}.
\]
Clearly,  $V_1$ is an upper bound of the latest birth in the first generation;
$V_1+V_2$ is an upper bound of the latest birth until the second generation and by induction 
\begin{equation}\label{eq:V}
V:=\sum_{g\geq 1}V_g
\end{equation} 
is clearly an upper bound for $L$. 
Hereon, we suppose
\begin{equation}\label{hp:densityCramer}
\text{$\mathbb{E}[B_{1,1,1}^k]<\infty$,  for any $k\in\mathbb N$,}
\end{equation}
and denote by $G_{P_{1,1}}$ the probability generating function of $P_{1,1}$. In the next proposition, we provide the fluctuations and the sharp deviations
of the number of points $S_t$ up to time $t$. 

\begin{Proposition}\label{prop:GW}
Assume \eqref{eq:subcritico}, \eqref{hp:thetac} and \eqref{hp:densityCramer}.
Then, setting $H(s,M_1):=N_1([0,s])$, $s\in [0,t]$, we have that:\\
\noindent$(i)$ The formulas \eqref{eq:CLT}, \eqref{eq:extCLT} and
\eqref{eq:fluctuationderivatives} of Theorem \ref{thm:fluctuations} hold. \\
\noindent$(ii)$ For any arbitrarily fixed $\sigma\in\mathbb N$, the formulas \eqref{eq:deviationordertlattice} and \eqref{eq:deviationordertlatticetail} 
of Theorem \ref{thm:deviationslattice}  hold.\\
\noindent$(iii)$ Setting 
\[
a_c:=\sup_{\theta\geq 0}(\theta-\log\mathbb{E}[\mathrm{e}^{\theta P_{1,1}}]),
\]
we have $a_c\in (0,a_Z]$.\\
\noindent$(iv)$ Set $b_c:=\sup\Theta$, where
\[
\Theta:=\{\theta:\,\,\theta<a_c\text{ and }\mathbb{E}[\mathrm{e}^{\theta Z}]<\mathrm{e}^b\},
\]
and $b$ is defined by \eqref{hp:thetac} $($note that $b_c>0$$)$.  We have
\[
\mathbb{E}[\mathrm{e}^{\theta Z}]=\mathrm{e}^{\theta}G_{P_{1,1}}(\mathbb{E}[\mathrm{e}^{\theta Z}])<\infty,\quad\text{for any $\theta\in (-\infty,b_c)$,}
\]
and the moments of $Z$ satisfy the formula: 
\begin{align}
\mathbb{E}[Z^n]&=1+
\sum_{k=1}^{n}k!\binom{n}{k}
\sum_{i=1}^{k}\frac{\mathbb{E}[P_{1,1}(P_{1,1}-1)\ldots(P_{1,1}-(i-1))]}{i!}\nonumber\\
&\qquad\qquad\qquad
\sum_{m_1+m_2+\ldots+m_i=k}
\frac{\mathbb{E}[Z^{m_1}]}{m_1!}\ldots\frac{\mathbb{E}[Z^{m_i}]}{m_i!},\quad\text{for any $n\geq 1$,}\label{eq:recursive2901}
\end{align}
where the third sum is taken over all the $m_1,\ldots,m_i\in\mathbb N$ such that $m_1+\ldots+m_i=k$.\\
\noindent$(v)$ For any $x\in(0,\lambda\mathbb{E}[Z\mathrm{e}^{b_c Z}])$, let
$\varrho_x\in (0,\mathbb{E}[\mathrm{e}^{b_c Z}])$ be the unique solution of the equation in $\varrho$:
\begin{equation}\label{eq:equationSab}
\frac{\lambda\varrho}{1-\varrho\frac{G'_{P_{11}}(\varrho)}{G_{P_{1,1}}(\varrho)}}=x.
\end{equation}
Then 
\[
\theta_x=\log\frac{\varrho_x}{G_{P_{1,1}}(\varrho_x)},\quad\varrho_x=\mathbb{E}[\mathrm{e}^{\theta_x Z}],
\]
and
\[
\mathbb{E}[Z^2\mathrm{e}^{\theta_x Z}]=
\varrho_x
\frac{1+\mathrm{e}^{\theta_x}(G'_{P_{1,1}}(\varrho_x)+
G''_{P_{1,1}}(\varrho_x))}
{(1-\mathrm{e}^{\theta_x}G'_{P_{1,1}}(\varrho_x))^2}.
\]
\noindent$(vi)$ The following estimates hold:\\
\noindent$(1)$ For any $x\in(0,\lambda(1-\mathbb{E}[P_{1,1}])^{-1})$ and any $q,q'>1$ such that $q^{-1}+q'^{-1}=1$, we have
\begin{equation}\label{eq:lowerbd2800BW}
p_1\mathbb{E}[B_{1,1,1}]
\mathrm{e}^{\theta_x}(1-\mathbb{E}[\mathrm{e}^{\theta_x Z}])\leq\varphi(\theta_x)
\leq-\theta_x\|B_{1,1,1}\|_{L^{q}(\mathbb P)}\|Z\|_{L^{q}(\mathbb P)}\|Z\|_{L^{q'}(\mathbb P)}.
\end{equation}
Note that these bounds are non trivial for any choice of the conjugate exponents.\\
\noindent$(2)$ For any $x\in(\lambda(1-\mathbb{E}[P_{1,1}])^{-1},\lambda\mathbb{E}[Z\mathrm{e}^{b_c Z}])$ and any $q,q',q_1,q_2>1$ such that
$q^{-1}+q'^{-1}=1$ and $q_1^{-1}+q_2^{-1}=1$, we have
\begin{align}
-\theta_x\|\mathrm{e}^{\theta_x Z}\|_{L^q(\mathbb P)}\|B_{1,1,1}\|_{L^{q'q_1}(\mathbb P)}\|Z\|_{L^{q'q_1}(\mathbb P)}
\|Z\|_{L^{q'q_2}(\mathbb P)}
&\leq\varphi(\theta_x)
\leq p_1\mathbb{E}[B_{1,1,1}]\mathrm{e}^{\theta_x}(1-\mathbb{E}[\mathrm{e}^{\theta_x Z}]).\label{eq:lowerbd2801BW}
\end{align}
Note that the upper bound is always non trivial. The lower bound is non trivial for any $q\in(1,b_c/\theta_x)$.

Here,  the function $\varphi(\cdot)$ involved in the relations \eqref{eq:lowerbd2800BW} and  \eqref{eq:lowerbd2801BW} 
is defined as in the statement of Proposition \ref{prop:PoissonCluster}. 
\end{Proposition}

\subsection{Example 1: Binomial offspring distribution}

Suppose that $P_{1,1}$ has a binomial distribution with parameters $(m,p)$ such that $\mathbb{E}[P_{1,1}]=mp<1$, $m\geq 1$, $p\in (0,1)$. Then condition
\eqref{eq:subcritico} is satisfied, and the assumption \eqref{hp:thetac} holds with $b=+\infty$, indeed
\[
G_{P_{1,1}}(\varrho)=(\varrho p+(1-p))^m,\quad\varrho\in\mathbb R.
\]
A standard computation gives
\[
a_c=b_c=\log\left(\frac{1}{mp}\left(\frac{m-1}{m(1-p)}\right)^{m-1}\right)
\]
(since $0^0:=1$, if $m=1$ then $a_c=b_c=-\log p$).
Using the formula \eqref{eq:totoff}, we have that $Z$ is distributed according to the Consul distribution, i.e.,
\[
\mathbb{P}(Z=k)=\frac{1}{k}\binom{km}{k-1}p^{k-1}(1-p)^{km-k+1},\quad k\in\mathbb N
\]
(see e.g. \cite{ConsulFamoye}). Therefore
\begin{align}
\mathbb{E}[\mathrm{e}^{b_c Z}]&=\sum_{k\geq 1}
\left(\frac{1}{mp}\left(\frac{m-1}{m(1-p)}\right)^{m-1}\right)^k
\frac{1}{k}\binom{km}{k-1}p^{k-1}(1-p)^{km-k+1}\nonumber\\
&=\frac{1-p}{p}\sum_{k\geq 1}\frac{1}{k}\binom{mk}{k-1}\left(\frac{(m-1)^{m-1}}{m^m}\right)^k=+\infty,\label{eq:diverge2901}
\end{align}
and so $\mathbb{E}[Z\mathrm{e}^{b_c Z}]=+\infty$. To check that the infinite sum in \eqref{eq:diverge2901} diverges, we note that it is trivially equal to $+\infty$
if $m=1$ and that, by construction, the total progeny of a Galton-Watson process with offspring distribution the binomial law with parameters $(m,p)$
is bigger than or equal to the total progeny of a Galton-Watson process with offspring distribution the Bernoulli law with parameter $p$.
Note that the equation \eqref{eq:equationSab} reads
\[
\frac{\lambda p \varrho^2+\lambda(1-p)\varrho}{p(1-m)\varrho+(1-p)}=x,
\]
i.e.,
\[
\lambda p \varrho^2+(\lambda(1-p)+(m-1)px)\varrho-(1-p)x=0,\quad x\in(0,\infty),
\]
which gives
\[
\varrho_x=\frac{-(\lambda(1-p)+(m-1)p x)+\sqrt{(\lambda(1-p)+(m-1)px)^2+4xp(1-p)\lambda}}{2\lambda p}.
\]
Therefore, we know explicitly the quantities $\theta_x$, $\mathbb{E}[\mathrm{e}^{\theta_x Z}]$ and $\mathbb{E}[Z^2\mathrm{e}^{\theta_x Z}]$.

Note that the formulas \eqref{eq:latticeugbreve} and \eqref{eq:latticemaggiugbreve} involve also the integral $\varphi(\theta_x)$
whose estimates \eqref{eq:lowerbd2800BW} and \eqref{eq:lowerbd2801BW} involve, in turn, the moments of $Z$ and of
$\mathrm{e}^{\theta_x Z}$. Hereafter, for the sake of completeness we briefly describe how those estimates can be computed.
For any $x\in(0,\lambda(1-mp)^{-1})$, choose, for instance,  $q=q'=2$
in the upper bound of \eqref{eq:lowerbd2800BW}. Then both the lower and the upper bounds on the integral are explicit (once the law of $B_{1,1,1}$ is fixed). 
Indeed, $p_1=mp(1-p)^{m-1}$, we already computed the quantities $\theta_x$ and $\mathbb{E}[\mathrm{e}^{\theta_x Z}]$ and by
\eqref{eq:recursive2901} one has $\mathbb{E}[Z]=(1-\mathbb{E}[P_{1,1}])^{-1}$, and
\begin{align}
\mathbb{E}[Z^2]&=1+2\mathbb{E}[P_{1,1}]\mathbb{E}[Z]+\mathbb{E}[P_{1,1}]\mathbb{E}[Z^2]+(\mathbb{E}[P_{1,1}^2]-\mathbb{E}[P_{1,1}])(\mathbb{E}[Z])^2,\nonumber
\end{align}
from which $\mathbb{E}[Z^2]$ is readily calculated.
For any $x\in(\lambda(1-mp)^{-1},\infty)$,
we note that  
(once the law of $B_{1,1,1}$ is fixed) the upper bound in \eqref{eq:lowerbd2801BW} is always explicit. As far as the lower bound  in \eqref{eq:lowerbd2801BW} is concerned,  we note that,  fixed $q\in (1,b_c/\theta_x)$, the finite quantity $\|\mathrm{e}^{\theta_x Z}\|_{L^q(\mathbb P)}$ can be computed e.g.  numerically (the law of $Z$ is 
explicitely known). Taking e.g. $q_1=q_2=2$ and letting $q'$ be the conjugate exponent of $q$, we have
\[
\|Z\|_{2q'}\leq\mathbb{E}[Z^{\lceil 2q'\rceil}]^{1/(2q')}
\]
(where $\lceil\cdot\rceil$ denotes the ceiling function),
which yields a computable lower bound (recall that $-\theta_x<0$), using the formula \eqref{eq:recursive2901}.

\subsection{Example 2: Geometric offspring distribution}

Suppose that $P_{1,1}$ has a geometric distribution with parameter $p\in (1/2,1)$. Then condition
\eqref{eq:subcritico} is satisfied,  and the assumption \eqref{hp:thetac} holds with $b=-\log(1-p)$, indeed
\[
G_{P_{1,1}}(\varrho)=\frac{p}{1-\varrho(1-p)},\quad |\varrho|<(1-p)^{-1}.
\]
A standard computation gives
\[
a_c=-\log (4p(1-p)).
\]
Using the formula \eqref{eq:totoff}, we have
\[
\mathbb{P}(Z=k)=\frac{1}{k}\binom{2(k-1)}{k-1}(1-p)^{k-1}p^k,\quad k\in\mathbb N.
\]
For all $\theta\leq a_c=-\log (4p(1-p))$,  it holds
\begin{align}
\sum_{k\geq 1}\frac{1}{k}\binom{2(k-1)}{k-1}[(1-p)p\mathrm{e}^{\theta}]^k&=
(1-p)p\mathrm{e}^{\theta}\sum_{h\geq 0}\frac{1}{h+1}\binom{2h}{h}[(1-p)p\mathrm{e}^{\theta}]^h\nonumber\\
&=\frac{1-\sqrt{1-4(1-p)p\mathrm{e}^{\theta}}}{2}<1,\nonumber
\end{align}
where the latter equality follows reconizing the generating function of the Catalan numbers.
Therefore,  for all $\theta\leq a_c$ we have $\mathbb{E}[\mathrm{e}^{\theta Z}]<1/(1-p)$,  which implies $b_c=a_c$.
Moreover,
\begin{align}
\mathbb{E}[\mathrm{e}^{b_c Z}]&=\frac{1}{4(1-p)}\sum_{k\geq 1}\frac{1}{k}\binom{2(k-1)}{k-1}\left(\frac{1}{4}\right)^{k-1}=\frac{1}{8(1-p)},\nonumber
\end{align}
and
\[
\mathbb{E}[Z\mathrm{e}^{b_c Z}]=\frac{1}{4(1-p)}\sum_{h\geq 0}\binom{2h}{h}\left(\frac{1}{4}\right)^{h}=+\infty.
\]
Note that the equation \eqref{eq:equationSab} reads
\[
\frac{\lambda\varrho(1-(1-p)\varrho)}{1-2(1-p)\varrho}=x,
\]
i.e.,
\[
-\lambda(1-p)\varrho^2+(2(1-p)x+\lambda)\varrho-x=0,\quad x>0,
\]
which gives
\[
\varrho_x=\frac{(2(1-p)x+\lambda)-\sqrt{(2(1-p)x+\lambda)^2-4\lambda x(1-p)}}{2\lambda(1-p)}.
\]
Therefore, we know explicitly the quantities $\theta_x$, $\mathbb{E}[\mathrm{e}^{\theta_x Z}]$ and $\mathbb{E}[Z^2\mathrm{e}^{\theta_x Z}]$.
Similar considerations as in the case of a binomial offspring distribution yield estimates of the integral $\varphi(\theta_x)$, which also appears
in the formulas \eqref{eq:latticeugbreve} and \eqref{eq:latticemaggiugbreve}.

\subsection{Example 3: Poisson offspring distribution, i.e., linear Hawkes processes with a general displacement distribution}

As already mentioned, precise deviations and fluctuations of classical Hawkes processes (or linear Hawkes processes) have been  studied in \cite{GZ}. Here, for the sake of completeness, we 
briefly explain how to apply Proposition \ref{prop:GW} even to the linear Hawkes process.

We emphasize that to apply Proposition \ref{prop:GW} we need to assume that the displacement distribution, i.e., the law of $B_{1,1,1}$, has all  moments finite.
The results in \cite{GZ}, instead, require only the existence of some moments for the displacement distribution, which, however, 
must have a density with respect to the Lebesgue measure. 
Hereon, instead,  we do not require the law of $B_{1,1,1}$ to be  absolutely continuous with respect to the Lebesgue measure.

So, suppose that $P_{1,1}$ has a Poisson distribution with mean $\mu\in (0,1)$. Then condition
\eqref{eq:subcritico} is satisfied,  and the assumption \eqref{hp:thetac} holds with $b=+\infty$, indeed
\[
G_{P_{1,1}}(\varrho)=\mathrm{e}^{\mu(\varrho-1)},\quad y\in\mathbb{R}.
\]
A standard computation gives
\[
a_c=\mu-1-\log\mu.
\]
Using the formula \eqref{eq:totoff}, we have that $Z$ is distributed according to the Borel distribution, i.e.,
\[
\mathbb{P}(Z=k)=\frac{(\mu k)^{k-1}}{k!}\mathrm{e}^{-\mu k},\quad k\in\mathbb N
\]
(see e.g.  \cite{ConsulFamoye}). Clearly, $b_c=a_c$. Moreover, standard computations yield
\begin{align}
\mathbb{E}[\mathrm{e}^{b_c Z}]=\mu^{-1}\quad\text{and}\quad\mathbb{E}[Z\mathrm{e}^{b_c Z}]=+\infty.\nonumber
\end{align}
Note that the equation \eqref{eq:equationSab} reads
\[
\frac{\lambda\varrho}{1-\mu\varrho}=x,
\]
which gives
\[
\varrho_x=\frac{x}{\lambda+\mu x}.
\]
Therefore, we know explicitly the quantities $\theta_x$, $\mathbb{E}[\mathrm{e}^{\theta_x Z}]$ and $\mathbb{E}[Z^2\mathrm{e}^{\theta_x Z}]$.
Similar considerations as in the case of a binomial offspring distribution yield estimates of the integral $\varphi(\theta_x)$, which also appears
in the formulas \eqref{eq:latticeugbreve} and \eqref{eq:latticemaggiugbreve}.

\section{Further applications}\label{sec:insurance}

\subsection{Ruin probabilities of risk processes with delayed claims}

Consider an insurance company with initial capital $u>0$ and premium rate $c>0$. As already mentioned in the Introduction, the total claim amount, up to time $t>0$, due to Incurred But Not Reported Claims is often modeled by a Poisson shot noise process $S_t$ of the form \eqref{eq:PSN},
see \cite{KM0, KM}.  We recall that, in this context, it is assumed
\begin{equation}\label{eq:HINSURANCE}
\text{For any $m\in\mathrm{M}$, the function $H(\cdot,m)$ is non-negative and non-decreasing. }
\end{equation}
The ruin probability of the insurance company is clearly given by
\[
\psi_{IBNR}(u):=\mathbb{P}\left(\sup_{t>0}(S_t-ct)\geq u\right),\quad u>0.
\]
Letting
\[
\psi_{CL}(u):=\mathbb{P}\left(\sup_{t>0}(C_t-ct)\geq u\right),\quad u>0,
\]
denote the ruin probability of the \lq\lq associated" Cram\'er-Lundberg risk process, where
\begin{equation*}
C_t:=\sum_{n\geq 1}H(\infty,M_n)\bold{1}_{(0,t]}(T_n),\quad t>0,
\end{equation*}
it is well-known that (see e.g. Theorem 1.2.2 in \cite{EKM}) if 
\begin{equation}\label{eqHinfinito}
\text{\eqref{eq:compLT} holds with $Z:=H(\infty,M_1)$ in place of $X$, $\lambda(\mathbb{E}[\mathrm{e}^{a Z}]-1)-ca>0$,} 
\end{equation}
\begin{equation}\label{eq:intxexpfinito}
\int_{0}^{\infty}x\mathrm{e}^{w x}\mathbb{P}(Z>x)\,\mathrm{d}x<\infty,
\end{equation}
and the \lq\lq net profit" condition 
\begin{equation}\label{eq:netprofit}
c>\lambda\mathbb{E}[Z]
\end{equation}
is satisfied, then
\begin{equation}\label{eq:CLestimate}
\psi_{CL}(u)=\left(\frac{\lambda w}{c-\lambda\mathbb{E}[Z]}\int_{0}^{\infty}x\mathrm{e}^{w x}\mathbb{P}(Z>x)\,\mathrm{d}x\right)^{-1}\mathrm{e}^{-w u}(1+o(1)),
\quad\text{as $u\to\infty$.}
\end{equation}
Here $w>0$ is the unique positive solution of the equation in $\gamma$:
\begin{equation}\label{eq:w}
\lambda(\mathbb{E}[\mathrm{e}^{\gamma Z}]-1)-c\gamma=0.
\end{equation}
Using large deviations techniques, it was proved in \cite{B} that, under the foregoing assumptions,
\[
\lim_{u\to\infty}\frac{1}{u}\log\psi_{IBNR}(u)=-w.
\]
The next proposition improves such ruin probability estimate.

\begin{Proposition}\label{cor:improvement}
Assume \eqref{eq:HINSURANCE},  \eqref{eqHinfinito}, \eqref{eq:intxexpfinito} and \eqref{eq:netprofit}. Then:\\
$(i)$ If \eqref{eq:nonlattice} and \eqref{eq:Lq} hold,  then
\begin{align}
&
\frac{\exp\left(\lambda\varphi(w)
\right)}
{w\sqrt{2\lambda\pi(\lambda\mathbb{E}[Z\mathrm{e}^{wZ}]-c)^{-1}\mathbb{E}[Z^2\mathrm{e}^{wZ}]}}\frac{\mathrm{e}^{-wu}}{\sqrt u}(1+o(1))
\leq\psi_{IBNR}(u)\nonumber\\
&\qquad\qquad\qquad
\leq
\left(\frac{\lambda w}{c-\lambda\mathbb{E}[Z]}\int_{0}^{\infty}x\mathrm{e}^{w x}\mathbb{P}(Z>x)\,\mathrm{d}x\right)^{-1}\mathrm{e}^{-w u}(1+o(1)),
\quad\text{as $u\to\infty$.}
\nonumber
\end{align}
$(ii)$ If \eqref{eq:integerpos} and \eqref{eq:Lqspeed} hold, then
\begin{align}
&\frac{\exp\left(\lambda\varphi(w)
\right)}
{w(1-\mathrm{e}^{-w})\sqrt{2\lambda\pi(\lambda\mathbb{E}[Z\mathrm{e}^{wZ}]-c)^{-1}\mathbb{E}[Z^2\mathrm{e}^{wZ}]}}\frac{\mathrm{e}^{-wu}}{\sqrt u}(1+o(1))
\leq\psi_{IBNR}(u)\nonumber\\
&\qquad\qquad\qquad
\leq
\left(\frac{\lambda w}{c-\lambda\mathbb{E}[Z]}\int_{0}^{\infty}x\mathrm{e}^{w x}\mathbb{P}(Z>x)\,\mathrm{d}x\right)^{-1}\mathrm{e}^{-w u}(1+o(1)),
\quad\text{as $u\to\infty$.}
\nonumber
\end{align}
\end{Proposition}

\subsubsection{Example 4: a risk process with delayed claims}\label{Ex:4}

We consider a risk process with delay in claim settlement and total claim amount at time $t>0$ given by the Poisson shot noise process
\[
S_t:=\sum_{n\geq 1}F(t-T_n)M_n\mathbf{1}_{(0,t]}(T_n),
\]
where $F$ is the distribution function of a law on $(0,\infty)$ with a finite mean. We suppose that the random variable $M_1$ is non-negative, it has a density with respect to the Lebesgue measure and it is such that $a_{M_1}\in (0,\infty]$ (see \eqref{eq:compLT}). Then, setting $\overline F:=1-F$,  we have
\[
M_1\int_0^\infty\overline{F}(s)\,\mathrm{d}s\in L^q(\mathbb P),\quad\text{for any $q>1$,}
\]
and the compound Poisson law $\phi_{\lambda,M_1}$ has a density with respect to the Lebesgue measure.  
Note that $Z$ is distributed as $M_1$ and the assumptions of Theorems 
\ref{thm:deviationsnonlattice} and \ref{thm:fluctuations} are satisfied.
\[
{\it Exponentially\,\,distributed\,\,claims:\,\,sharp\,\,deviations\,\,at\,\,scales\,\,O(t)}
\]
Suppose that $M_1$ is exponentially distributed with mean $1/a$, $a\in (0,\infty)$. Then,  a straightforward computation shows
\begin{equation}\label{eq:tetax}
\theta_x=a-\sqrt{\frac{\lambda a}{x}},\quad x>\lambda/a,
\end{equation}
\begin{equation}\label{eq:24DIC1}
\mathbb{E}[\mathrm{e}^{\theta_x M_1}]=\frac{a}{a-\theta_x}=\sqrt{\frac{ax}{\lambda}},
\end{equation}
\begin{equation}\label{eq:24DIC2}
\mathbb{E}[M_1^2 \mathrm{e}^{\theta_x M_1}]=\frac{2a}{(a-\theta_x)^3}=2\frac{x}{\lambda}
\sqrt{\frac{x}{\lambda a}},
\end{equation}
and
\begin{align}
\varphi(\theta_x)=\int_0^\infty(\mathbb{E}[\mathrm{e}^{\theta_x F(s)M_1}]-\mathbb{E}[\mathrm{e}^{\theta_x M_1}])\,\mathrm{d}s&=
\int_0^\infty\left(\frac{a}{a-\theta_x F(s)}-\frac{a}{a-\theta_x}\right)\,\mathrm{d}s\nonumber\\
&=-a\left(\sqrt{\frac{ax}{\lambda}}-1\right)\Phi(\lambda,a,x,F),\nonumber
\end{align}
where
\[
\Phi(\lambda,a,x,F):=\int_0^\infty\frac{\overline{F}(s)}{a-\left(a-\sqrt{\frac{\lambda a}{x}}\right)F(s)}\,\mathrm{d}s.
\]
So by Theorem \ref{thm:deviationsnonlattice} we have that, for any $x>\lambda/a$, as $t\to\infty$,
\[
\mathbb{P}(S_t\geq tx)=\frac{(\lambda/ax)^{1/4}}{2\sqrt{\pi}(\sqrt{ax}-\sqrt{\lambda})}\exp\left(-a(\sqrt{\lambda a x}-\lambda)\Phi(\lambda,a,x,F)\right)\frac{1}{\sqrt t}
\mathrm{e}^{-(xa+\lambda-2\sqrt{\lambda a x})t}(1+o(1)).
\]
For specific choices of $F(\cdot)$ the quantity $\Phi$ can be computed.
For instance, if $F(\cdot)$ is the distribution function of the uniform law on $(0,1)$, then
\[
\int_0^\infty\frac{\overline{F}(s)}{a-\left(a-\sqrt{\frac{\lambda a}{x}}\right)F(s)}\,\mathrm{d}s=-\frac{\sqrt{\lambda x/a}}
{(\sqrt{ax}-\sqrt{\lambda})^2}\log\sqrt{\frac{ax}{\lambda}}
+\frac{1}{a-\sqrt{\frac{\lambda a}{x}}}.
\]
For this model, one may easily compute the ruin probabilities estimates provided by Proposition \ref{cor:improvement}$(i)$. We omit the details.
\[
{\it Fluctuations\,\,at\,\,scales\,\,o(t)}
\]
Hereon,  by applying Theorem \ref{thm:fluctuations}$(iii)$ with $y(t):=o(t^{1/4})$,  $y(t)\to\infty$,
we study the fluctuations of $S_t$ around its asymptotic mean at a scale of order $o(t^{3/4})$. 
To perform the computation, we do not need to know the distribution of $M_1$
(as it happens to study the sharp deviations of $S_t$ from its asymptotic mean, where we need to compute the function $\theta_{\cdot}$), but only the first four moments of $M_1$.
Indeed,  under the foregoing assumptions, the formula \eqref{eq:fluctuationderivatives} holds with $Z=M_1$ and $m=4$ and, assuming
the knowledge of the quantities $\mathbb{E}[M_1^i]$, $i=1,2,3,4$, we only need to compute $\theta_{\lambda\mathbb{E}[M_1]}^{(j)}$ for $j=2,3$. To this aim,
by the recursive formula \eqref{eq:recursive} we have
\begin{equation*}
\theta_{\lambda\mathbb{E}[M_1]}^{(2)}=g'(0)\theta_{\lambda\mathbb{E}[M_1]}^{(1)}=-\frac{\lambda\mathbb{E}[M_1^3]}{(\lambda\mathbb{E}[M_1^2])^3}
\end{equation*}
and
\begin{align}
\theta_{\lambda\mathbb{E}[M_1]}^{(3)}&=g'(0)\theta_{\lambda\mathbb{E}[M_1]}^{(2)}+g''(0)(\theta_{\lambda\mathbb{E}[M_1]}^{(1)})^2\nonumber\\
&=\frac{(\lambda\mathbb{E}[M_1^3])^2}{(\lambda\mathbb{E}[M_1^2])^5}-\frac{\lambda}{(\lambda\mathbb{E}[M_1^2])^2}\left(\mathbb{E}[M_1^4]
-2\lambda\frac{(\mathbb{E}[M_1^3])^2}{\mathbb{E}[M_1^2]}\right).\nonumber
\end{align}

\subsection{A teletraffic model}

In this subsection we briefly discuss the application of Theorems \ref{thm:deviationsnonlattice} and \ref{thm:fluctuations}  to a teletraffic model proposed in \cite{KL}.  

We consider an infinite servers queuing system, handling in parallel active jobs
(connections).
New jobs arrive at the system according to a homogeneous Poisson process  
$\{T_n\}_{n\geq 1}$. 
Let $\{M_n\}_{n\geq 1}$  denote the sequence of the processing times of the jobs, which are assumed to be independent and identically distributed. Every  job en-queued is served by the system at 
unit rate.  At time $t$, the number of active jobs in the system is given by
\[
X_t:=\sum_{n\geq 1}\mathbf{1}_{(0,M_n]}(t-T_n)\mathbf{1}_{(0,t]}(T_n).
\]
Then,  the total workload processed by the system up to time $t$ is given by the Poisson shot noise process
\[
S_t:=\int_{0}^{t}X_s\,\mathrm{d}s=\sum_{n\geq 1}[(t-T_n)\wedge M_n]\bold{1}_{(0,t]}(T_n),
\]
where, for $a,b\in\mathbb{R}$, the symbol $a\wedge b$ denotes the minimum between $a$ and $b$.  
Hereon, we suppose that the random variable $M_1$ has a density with respect to the Lebesgue measure and 
that $a_{M_1}\in (0,\infty]$ (see \eqref{eq:compLT}). Then
\[
\int_0^\infty(M_1-s\wedge M_1)\,\mathrm{d}s=\frac{M_1^2}{2}\in L^q(\mathbb P),\quad\text{for any $q>1$,}
\]
and the compound Poisson law $\phi_{\lambda,M_1}$ is absolutely continuous with respect to the Lebesgue measure. 
Here again,  $Z$ is distributed as $M_1$ and the assumptions of Theorems 
\ref{thm:deviationsnonlattice} and \ref{thm:fluctuations} are satisfied.
Assuming the knowledge of only the first four moments of $M_1$, by the same application of Theorem \ref{thm:fluctuations}$(iii)$ with $y(t):=o(t^{1/4})$,
$y(t)\to\infty$, discussed in the Example 4, we can quantify the fluctuations of $S_t$ around its asymptotic mean at a scale of order $o(t^{3/4})$. 
If we assume that $M_1$ is exponentially distributed with mean $1/a$, $a\in (0,\infty)$, then the function $\theta_\cdot$ is given by
\eqref{eq:tetax} and the quantities $\mathbb{E}[\mathrm{e}^{\theta_x M_1}]$ and $\mathbb{E}[M_1^2 \mathrm{e}^{\theta_x M_1}]$ are provided by
\eqref{eq:24DIC1} and \eqref{eq:24DIC2}, respectively. Moreover
\begin{align}
\varphi(\theta_x)=\int_0^\infty(\mathbb{E}[\mathrm{e}^{\theta_x(s\wedge M_1)}]-\mathbb{E}[\mathrm{e}^{\theta_x M_1}])\,\mathrm{d}s&=\mathbb{E}\left[\int_0^{M_1}(\mathrm{e}^{\theta_x s}-\mathrm{e}^{\theta_xM_1})\,\mathrm{d}s\right]
\nonumber\\
&=-\frac{1}{\theta_x}+\frac{a}{\theta_x(a-\theta_x)}-\frac{a}{(a-\theta_x)^2}\nonumber\\
&=\sqrt{\frac{x}{\lambda a}}-\frac{x}{\lambda}.\nonumber
\end{align}
So by Theorem \ref{thm:deviationsnonlattice} we have the explicit asymptotic expression of the tail of the processed workload.
 
\section{Stable probability approximation of Poisson shot noise processes}\label{sec:stablegrosso}

We recall that throughout this paper we denote by $\phi_{c,\alpha,\beta}$ the stable law with scale parameter $c>0$, stability parameter $\alpha\in (0,2]$ and 
skewness parameter $\beta\in [-1,1]$ (see Subsection \ref{subsec:stable} for the definition of stable distribution and its first properties). 
In the context of the stable probability approximation of the Poisson shot noise, we consider the following assumptions on the model \eqref{eq:PSN}:
\begin{equation}\label{eq:stablemark}
\text{The random variable $M_1$ is distributed according to $\phi_{c,\alpha,\beta}$ with
either $\alpha\neq 1$ or $\alpha=1$ and $\beta=0$.}
\end{equation}
\begin{equation}\label{eq:multiplicative}
\text{For any $t>0$ and $m\in\mathbb{R}$, $H(t,m):=m F(t)$, for some distribution function $F:[0,\infty)\to [0,1]$.}
\end{equation}

The following theorem holds.

\begin{Theorem}\label{cor:coninlaw}
Assume \eqref{eq:stablemark} and \eqref{eq:multiplicative} and let $S$ be a random variable with law $\phi_{c\lambda^{1/\alpha},\alpha,\beta}$.
Then $S_t/t^{\frac{1}{\alpha}}$ converges in law to $S$, as $t\to\infty$. 
\end{Theorem}

The next theorem quantifies the speed of such weak convergence. 
For the sake of readability,  we state this result by giving estimates of the Kolmogorov distance
between $S_t/t^{1/\alpha}$ and $S$ in terms of big $O$ functions. 
However,  in the proof we compute explicitly the constants involved in the bounds; anyway,  such constants have a complicated expression.

Letting $X,Y$ denote two real-valued random variables,  we recall that the Kolmogorov distance between $X$ and $Y$
is defined by
\[
\mathrm{d}_{\mathrm{Kol}}(X,Y):=\sup_{x\in\mathbb R}|\mathbb{P}(X\leq x)-\mathbb{P}(Y\leq x)|.
\]

\begin{Theorem}\label{prop:speed}
Let the assumptions and notation of Theorem \ref{cor:coninlaw} prevail. Then:\\
\noindent $(i)$ Define
\begin{equation*}
A:=\Biggl\{\eta\in [\max\{0,(\alpha-1)/\alpha\},1/2]:\,\,\limsup_{t\to\infty}\frac{\int_0^t(1-F(s)^{\alpha})\,\mathrm{d}s}{t^\eta}<\infty\Biggr\}
\end{equation*}
and  $\eta_0:=\inf A$.\\
If
\begin{equation}\label{eq:sigma}
\eta_0\in A
\end{equation}
then 
\[
\mathrm{d}_{\mathrm{Kol}}(S_t/t^{1/\alpha},S)=O\left(\frac{1}{t^{1-\eta_0}}\right),\quad\text{as $t\to\infty$.}
\]
Instead, if
\[ 
A\neq\emptyset \quad \text{and} \quad  \eta_0\notin A,
\] then for any $\widetilde{\eta}_0\in (\eta_0,1)$ it holds
\[
\mathrm{d}_{\mathrm{Kol}}(S_t/t^{1/\alpha},S)=O\left(\frac{1}{t^{1-\widetilde{\eta}_0}}\right),\quad\text{as $t\to\infty$.}
\]
\noindent$(ii)$ If 
\begin{equation}\label{eq:sigma2}
\int_0^\infty(1-F(s)^{\alpha})\,\mathrm{d}s\in [0,\infty)
\end{equation}
then
\[
\mathrm{d}_{\mathrm{Kol}}(S_t/t^{1/\alpha},S)=O\left(\frac{1}{t^{1/\alpha}}\right)\bold{1}\{\alpha\in (1,2]\}+
O\left(\frac{1}{t}\right)\bold{1}\{\alpha\in (0,1)\cup\{\alpha=1,\,\beta=0\}\},\quad\text{as $t\to\infty$.}
\]
\end{Theorem}

\begin{Remark}
In \cite{KMS} the authors study the weak convergence to a multivariate $\alpha$-stable distribution of the finite-dimensional distributions
of a properly centred and scaled Poisson shot noise process. In particular, Section 3.3 of \cite{KMS} considers Poisson shot noise processes of the form
\[
X_t:=\sum_{n\geq 1}M_n h(t-T_n)\bold{1}_{(0,t]}(T_n),\quad t>0,
\]
where $\{T_n\}_{n\geq 1}$ is a homogeneous Poisson process on $(0,\infty)$, independent of the sequence of independent and identically distributed random variables $\{M_n\}_{n\geq 1}$, and $h(\cdot)$ is a measurable and positive function. By the theory developed in \cite{KMS} it follows that if the law of $M_1$ is regularly varying at infinity with a constant sign, and the shot shape function $h(\cdot)$ is bounded on the compacts and regularly varying too, then $X_t$, properly centred and scaled, converges weakly to a stable law.  Although the Poisson shot noise model $S_t$ considered in Theorem
\ref{cor:coninlaw} has a multiplicative noise, we consider different hypotheses on the mark $M_1$ and on the shot shape function $F(\cdot)$. For instance, 
we do not assume that both the distribution of $M_1$ and the shot shape function are regularly varying. Furthermore,
Theorem \ref{prop:speed} 
provides the speed of the weak convergence to the stable law, which is far from the achievements in \cite{KMS}.
\end{Remark}  

\subsection{Example 5: shot shapes for the stable probability approximation}

Assuming \eqref{eq:stablemark} and \eqref{eq:multiplicative}, hereon we discuss some possible choices of the shot shape $F(\cdot)$ which allows us  to apply Theorem \ref{prop:speed}.  

We start noticing that if
the law with distribution function $F(\cdot)$ has a finite mean, then by Theorem \ref{prop:speed}$(ii)$ we have
\[
\mathrm{d}_{\mathrm{Kol}}(S_t/t^{1/\alpha},S)=O\left(\frac{1}{t^{1/\alpha}}\right)\bold{1}\{\alpha\in (1,2]\}+
O\left(\frac{1}{t}\right)\bold{1}\{\alpha\in (0,1)\cup\{\alpha=1,\,\beta=0\}\}.
\]
Indeed, since $\alpha\in (0,2]$ we have $F(t)^2\leq F(t)^{\alpha}$ for any $t>0$. Therefore
\begin{align}
\int_0^\infty(1-F(s)^\alpha)\,\mathrm{d}s&\leq\int_0^\infty(1-F(s)^2)\,\mathrm{d}s
\leq2\int_0^\infty(1-F(s))\,\mathrm{d}s<\infty.\nonumber
\end{align}

If either $\alpha\in (1,2]$ or $\alpha=1$ and $\beta=0$, and the law with distribution function 
$F(\cdot)$ has an infinite mean,  then Theorem \ref{prop:speed}$(ii)$ can not be applied, indeed $F(t)\geq F(t)^\alpha$ for any $t>0$. In such a case we can hope to apply
Theorem \ref{prop:speed}$(i)$. This is the case, for instance,  when $F(\cdot)$ is the distribution function of a Pareto law with parameter $1$, i.e.,
\[
F(t):=\frac{t}{1+t}\bold{1}_{[0,\infty)}(t).
\]
Indeed, we have
\[
\int_0^\infty(1-F(t))\,\mathrm{d}t=\infty,
\]
and so Theorem \ref{prop:speed}$(ii)$ does not apply (condition \eqref{eq:sigma2} is not satisfied). However, by applying twice de l'Hopital's theorem, for any 
$\eta\in (0,1/2]$ we have
\begin{align}
\lim_{t\to\infty}\frac{\int_0^t(1-F(s)^\alpha)\,\mathrm{d}s}{t^\eta}&=\lim_{t\to\infty}\frac{1-F(t)^\alpha}{\eta t^{\eta-1}}\nonumber\\
&=\frac{\alpha}{\eta(1-\eta)}\lim_{t\to\infty}\frac{F(t)^{\alpha-1}\frac{t^2}{(t+1)^2}}{t^\eta}=0.\label{eq:pareto0}
\end{align}
If $\alpha\in (1,2]$, then 
$\eta_0=(\alpha-1)/\alpha\in A$ (the set $A$ is defined in the statement of Theorem \ref{prop:speed}) and so
by Theorem \ref{prop:speed}$(i)$ we have 
$d_{\mathrm{Kol}}(S_t/t^{1/\alpha},S)=O\left(1/t^{1/\alpha}\right)$. If $\alpha=1$ and $\beta=0$, then $\eta_0=0\notin A$ and so
by Theorem \ref{prop:speed}$(i)$ we have 
$d_{\mathrm{Kol}}(S_t/t,S)=O\left(1/t^{1-\eta}\right)$, for any $\eta\in (0,1)$.
Of interest for applications in insurance it is the case when the total pay-off $M_1$ has the heavy tail L\'evy distribution, which corresponds to the stable law
with parameters $c=1$, $\alpha=1/2$ and $\beta=1$. In such a case if the delay in claim settlement is modeled by a distribution function $F(\cdot)$ with a finite mean, then we fall in the first case discussed above and so
\[
\mathrm{d}_{\mathrm{Kol}}(S_t/t^2,S)=O\left(\frac{1}{t}\right).
\]
If, for instance, $F(\cdot)$ is the distribution function of a Pareto law with parameter $1$, then (by the some computation as in \eqref{eq:pareto0})
\[
\lim_{t\to\infty}\frac{\int_0^t(1-F(s)^{1/2})\,\mathrm{d}s}{t^\eta}=0,\quad\text{for any $\eta\in (0,1/2]$,}
\]
and so $\eta_0=0\notin A$. Therefore
by Theorem \ref{prop:speed}$(i)$ we have 
$d_{\mathrm{Kol}}(S_t/t^2,S)=O\left(1/t^{1-\eta}\right)$, for any $\eta\in (0,1)$.  

\section{Proofs}\label{sec:proofs}

\subsection{Proof of Theorem \ref{thm:deviationsnonlattice}}\label{sec:precise}

As already mentioned in the Introduction,  the proof is based on the mod-$\phi$ convergence theory. 
So, first we state a lemma which concerns the mod-compound Poisson convergence of Poisson shot noise processes (i.e., the mod-$\phi$ convergence of Poisson shot noise processes when $\phi$ is the compound Poisson law),  then
we prove Theorem \ref{thm:deviationsnonlattice} applying the results in \cite{FMN}.  

The following lemma, whose proof is postponed at the end of the subsection, holds.

\begin{Lemma}\label{thm:modphiR}
Assume \eqref{eq:lighttailCP}, \eqref{eq:noneg} and \eqref{eq:Lq}. Then $\{S_t\}_{t>0}$ converges mod-$\phi_{\lambda,Z}$ on $D_{cp}:=D_{cp}(a)$ 
with parameter function $x(t):=t$ and limiting function
\begin{equation}\label{eq:psi}
\psi(z):=\mathrm{e}^{\lambda\varphi(z)},
\end{equation}
where $\varphi(\cdot)$ is defined by \eqref{eq:varphi27Lug2021}.
\end{Lemma}

Before proving Theorem \ref{thm:deviationsnonlattice}, we note that, under the assumptions and notation
of Lemma \ref{thm:modphiR},  by the theory of mod-$\phi$ convergence (see p. 20 of \cite{FMN}) we have the asymptotic expressions of the mean and the variance of 
$S_t$ anticipated in
\eqref{eq:mediavarianza}. 
\\

\noindent{\it Proof\,\,of\,\,Theorem\,\,\ref{thm:deviationsnonlattice}.} Let $x\in(\lambda\mathbb{E}[Z],\lambda\mathbb{E}[Z\mathrm{e}^{a Z}])$ be fixed. Reasoning by contradiction, suppose 
\begin{equation}\label{eq:limsup7gen}
\limsup_{t\to\infty}\frac{\mathbb{P}(S_t\geq tx)}{
\frac{\exp\left(-t\eta_{\lambda,Z}^*(x)+\lambda\varphi(\theta_x)
\right)}{\theta_x\sqrt{2\lambda\pi t\mathbb{E}[Z^2\mathrm{e}^{\theta_x Z}]}}}\neq 1,
\end{equation}
and let $\{t_k\}_{k\geq 1}\subset (0,\infty)$ be a divergent sequence that realizes the $\limsup$. By Lemma \ref{thm:modphiR} the stochastic process
$\{S_{t_k}\}_{k\geq 1}$ converges mod-$\phi_{\lambda,Z}$ on $D_{cp}$ with parameter sequence $x(t_k):=t_k$ and limiting function
\eqref{eq:psi}. So, by Theorem 4.2.1 in \cite{FMN} easily follows that the relation \eqref{eq:deviationordert} holds with $t_k$ in place of $t$.
Indeed, in our framework, the functions $F(\cdot)$, $h_{\cdot}$ and $\eta(\cdot)$ in the Theorem 4.2.1 of \cite{FMN} are given, respectively, by
the functions $\eta_{\lambda,Z}^*(\cdot)$ (in \eqref{eq:FenchelLegendre}), $\theta_{\cdot}$ and $\eta_{\lambda,Z}(\cdot)$. Moreover, the function $\psi(\cdot)$ and the quantities $c$ and $d$ of Theorem 4.2.1 
in \cite{FMN} are given, respectively, by the function 
$\psi(\cdot)$ defined by \eqref{eq:psi}, and $c=-\infty$ and $d=a$.Therefore we reached a contradiction and in \eqref{eq:limsup7gen} we have the equality.
The same arguments hold if in \eqref{eq:limsup7gen} we replace the $\limsup$ with the $\liminf$.
The proof is completed.\\
\noindent$\square$

\noindent{\it Proof\,\,of\,\,Lemma \ref{thm:modphiR}.} Note that $\{(T_n,M_n)\}_{n\geq 1}$ is a Poisson process on $(0,\infty)\times\mathrm{M}$ with mean measure
$\lambda\mathrm{d}t\mathbb{P}_{M_1}(\mathrm{d}m)$. So, by the expression of the Laplace functional of a Poisson process (see e.g. \cite{DVJ}), we have
\[
\mathbb{E}[\mathrm{e}^{z S_t}]=\exp\left(\lambda\int_0^t(\mathbb{E}[\mathrm{e}^{zH(s,M_1)}]-1)\,\mathrm{d}s\right),\quad\text{for any $t>0$ and $z\in D_{cp}$.}
\]
By the elementary relation $|\mathrm{e}^{z}|=\mathrm{e}^{\mathrm{Re}z}$, $z\in\mathbb{C}$,
we have
\begin{align}
|\mathbb{E}[\mathrm{e}^{z S_t}]|&\leq\mathbb{E}[\mathrm{e}^{(\mathrm{Re}z) S_t}]\nonumber\\
&=\exp\left(\lambda\int_0^t(\mathbb{E}[\mathrm{e}^{(\mathrm{Re}z)H(s,M_1)}]-1)\bold{1}\{\mathrm{Re}z\leq 0\}\,\mathrm{d}s+
\lambda\int_0^t(\mathbb{E}[\mathrm{e}^{(\mathrm{Re}z)H(s,M_1)}]-1)\bold{1}\{\mathrm{Re}z>0\}\,\mathrm{d}s\right)\nonumber\\
&\leq\exp(\lambda t(\mathbb{E}[\mathrm{e}^{(\mathrm{Re}z)Z}]-1)\bold{1}\{\mathrm{Re}z>0\})<\infty,\quad\text{for any $t>0$ and $z\in D_{cp}$.}\nonumber
\end{align}
Therefore, the Laplace transform of $S_t$, $t>0$, is defined on $D_{cp}$. 
We now check that $\psi(\cdot)$, defined by \eqref{eq:psi}, does not vanish on $\mathrm{Re}D_{cp}=(-\infty,a)$ and it is analytic on $D_{cp}$. By the
Lemma \ref{le:standardineq}, for any $z\in D_{cp}$ and $s>0$, we have
\begin{align}
|\mathrm{e}^{z H(s,M_1)}-\mathrm{e}^{z Z}|&\leq |z|(Z-H(s,M_1))\mathrm{e}^{\max\{(\mathrm{Re}z)H(s,M_1),(\mathrm{Re}z)Z\}}\nonumber\\
&=|z|(Z-H(s,M_1))\mathrm{e}^{\max\{(\mathrm{Re}z)H(s,M_1),(\mathrm{Re}z)Z\}}\mathbf{1}\{\mathrm{Re}z\leq 0\}\nonumber\\
&\qquad\qquad
+|z|(Z-H(s,M_1))\mathrm{e}^{\max\{(\mathrm{Re}z)H(s,M_1),(\mathrm{Re}z)Z\}}\mathbf{1}\{\mathrm{Re}z>0\}\nonumber\\
&\leq 
|z|(Z-H(s,M_1))\mathrm{e}^{\max\{0, \mathrm{Re}z\}Z}
\quad\text{$\mathbb{P}$-a.s..}\nonumber
\end{align}
Therefore, for any $z\in D_{cp}$ and $s>0$,
\begin{align}
\mathbb{E}[|\mathrm{e}^{z H(s,M_1)}-\mathrm{e}^{z Z}|]&\leq
|z|\mathbb{E}[(Z-H(s,M_1))\mathrm{e}^{\max\{0,\mathrm{Re}z\}Z}].
\label{eq:moduliexp}
\end{align}
So, applying H\"older's inequality with 
$p\in(1,a/\max\{0,\mathrm{Re}z\})$,
where we conventionally set 
$a/0:=\infty$, and $q>1$ such that $p^{-1}+q^{-1}=1$
(note that the conjugate exponents depend on $z$), we have
\begin{align}
|\varphi(z)|
&\leq |z|
\Big\|\int_0^\infty(Z-H(s,M_1))\,\mathrm{d}s\Big\|_{L^q(\mathbb P)}\Big\|\mathrm{e}^{\max\{0, \mathrm{Re}z\}Z}\Big\|_{L^p(\mathbb P)}<\infty 
,\quad\text{for any $z\in D_{cp}$,}\nonumber
\end{align}
where the latter term is finite due to the assumptions \eqref{eq:lighttailCP} and \eqref{eq:Lq}. So 
\[
\psi(x)=\mathrm{e}^{\lambda\varphi(x)}>0,\quad\text{for any $x\in(-\infty,a)$.}
\]
To prove that $\psi(\cdot)$ is analytic on $D_{cp}$ we show that $\varphi(\cdot)$ is such, and to this aim we prove that
\begin{equation}\label{eq:varphintovarphi}
\sup_{z\in K}|\varphi_t(z)-\varphi(z)|\to 0,\quad\text{as $t\to\infty$, for any compact $K\subset D_{cp}$,}
\end{equation}
where
\[
\varphi_t(z):=\int_0^t(\mathbb{E}[\mathrm{e}^{z H(s,M_1)}]-\mathbb{E}[\mathrm{e}^{z Z}])\,\mathrm{d}s,\quad\text{$z\in D_{cp}$}
\]
(then the analiticity of $\varphi(\cdot)$ on $D_{cp}$ it follows by the analiticity of $\varphi_t$ on $D_{cp}$ for any $t>0$). 
Let $K\subset D_{cp}$ be an arbitrarily fixed compact. Without loss of generality, we suppose that there exists $a'\in (0,a)$ such that
$K\cap\{z:\,\,0<\mathrm{Re}z\leq a'\}\neq\emptyset$. Indeed, if $K\subset\{z:\,\,\mathrm{Re}z\leq 0\}$, then there exists $K'\subset D_{cp}$ compact such that
$K\subset K'$ and $K'\cap\{z:\,\,0<\mathrm{Re}z\leq a'\}\neq\emptyset$. Hereon, we set
\[
m_K:=\max(\mathrm{Re}K).
\]
By construction, $0<m_K<a$, and, setting $\overline{M}_K:=\sup_{z\in K}|z|$, by \eqref{eq:moduliexp} it follows
\begin{equation*}
\mathbb{E}[|\mathrm{e}^{z H(s,M_1)}-\mathrm{e}^{z Z}|]\leq\overline{M}_K\mathbb{E}\left[(Z-H(s,M_1))
\mathrm{e}^{m_K Z}\right],
\quad\text{$\forall$ $s>0$ and $z\in K$.}
\end{equation*}
By this relation and the definition of $\varphi_t$ and $\varphi$, for any $t>0$ and $z\in K$, we have
\begin{align}
|\varphi_t(z)-\varphi(z)|
&\leq\overline{M}_K\mathbb{E}\left[\mathrm{e}^{m_K Z}\int_t^\infty(Z-H(s,M_1))\,\mathrm{d}s
\right].\nonumber
\end{align}
Applying  H\"older's inequality with $p\in (1,a/m_K)$ and $q>1$ such that $p^{-1}+q^{-1}=1$
(note that in this case the conjugate exponents depend on $K$, but not on $z$) 
for any $t>0$ and $z\in K$, we have
\begin{align}
|\varphi_t(z)-\varphi(z)|
&\leq \mathcal{M}_K(t):=\overline{M}_K\left(\mathbb{E}[\mathrm{e}^{p\,m_K Z}]^{1/p}\right)\Big\|\int_t^\infty(Z-H(s,M_1))\,\mathrm{d}s\Big\|_{L^q(\mathbb P)}.
\label{eq:moduliexp1}
\end{align}
The claim \eqref{eq:varphintovarphi} follows taking first the supremum over $K$ on this latter relation, and then letting $t$ tend to infinity. Indeed,
by the assumption \eqref{eq:Lq} and the dominated convergence theorem we immediately have
\[
\Big\|\int_t^\infty(Z-H(s,M_1))\,\mathrm{d}s\Big\|_{L^q(\mathbb P)}^q\to 0,\quad\text{as $t\to\infty$.}
\]
Finally we prove that
\[
\psi_t(z):=\mathrm{e}^{-t\eta_{\lambda,Z}(z)}\mathbb{E}[\mathrm{e}^{z S_t}],\quad\text{$t>0$, $z\in D_{cp}$,} 
\]
converges to $\psi(\cdot)$ uniformly on the compacts of $D_{cp}$.
Let $K\subset D_{cp}$ be an arbitrarily fixed compact. For any $z\in K$, by first applying Lemma \ref{le:standardineq}, 
then using the elementary inequality $\mathrm{Re}w\leq |w|$, $w\in\mathbb{C}$, and finally exploiting
\eqref{eq:moduliexp1}, we have
\begin{align}
&|\psi_t(z)-\psi(z)|=\Big|\exp\left(\lambda\int_0^{t}(\mathbb{E}[\mathrm{e}^{z H(s,M_1)}]-\mathbb{E}[\mathrm{e}^{z Z}])\,\mathrm{d}s\right)
-\exp(\lambda\varphi(z))
\Big|\nonumber\\
&\qquad\leq\lambda\int_{t}^{\infty}\mathbb{E}[|\mathrm{e}^{z H(s,M_1)}-\mathrm{e}^{z Z}|]\,\mathrm{d}s\times
\exp\left(\lambda\int_0^{\infty}\mathbb{E}[|\mathrm{e}^{z H(s,M_1)}-\mathrm{e}^{z Z}|]\,\mathrm{d}s\right)
\leq\left(\lambda\mathrm{e}^{\lambda\mathcal{M}_K(0)}\right)\mathcal{M}_K(t).\label{eq:7gen1}
\end{align}
The claimed local uniform convergence follows by first taking the supremum over $K$ on this relation and then letting $t$ tend to infinity.
\\
\noindent$\square$ 

\subsection{Proof of Theorem \ref{thm:deviationslattice}}\label{sec:preciselattice}

Here again, the proof is based on the mod-$\phi$ convergence theory.  However, we need to refine Lemma 
\ref{thm:modphiR} since we need the mod-compound Poisson convergence
with a polynomial decay (see Definition \ref{def:mod}). So, first we state such a refinement in a lemma and then we prove Theorem \ref{thm:deviationslattice}
applying the results in \cite{FMN}.  

The following lemma, whose proof is postponed at the end of the subsection, holds.

\begin{Lemma}\label{thm:modphiRspeed}
Assume \eqref{eq:lighttailCP}, \eqref{eq:integerpos} and \eqref{eq:Lqspeed}. Then $\{S_t\}_{t>0}$ converges mod-$\phi_{\lambda,Z}$ on $D_{cp}:=D_{cp}(a)$ 
with speed $O(t^{-\sigma})$ and limiting function \eqref{eq:psi}.
\end{Lemma}

\noindent{\it Proof\,\,of\,\,Theorem\,\,\ref{thm:deviationslattice}.} 
Although the proof is conceptually similar to the proof of Theorem \ref{thm:deviationsnonlattice}, since it exploits different results of the mod-$\phi$
convergence theory, we provide some details.\\
\noindent{\it Proof\,\,of\,\,Part\,\,$(i)$.}\\
Let $x\in (0,\lambda\mathbb{E}[Z\mathrm{e}^{a Z}])$ be such that $tx\in\mathbb{N}$. Reasoning by contradiction, suppose 
\begin{equation}\label{eq:limsup7genSBR}
\limsup_{t\to\infty}t^{-\sigma}\frac{\Big|\mathbb{P}(S_t=tx)-
\frac{\exp(-t\eta_{\lambda,Z}^*(x))}{\sqrt{2\lambda\pi t\mathbb{E}[Z^2\mathrm{e}^{\theta_x Z}]}}
\left(\psi(\theta_x)+\sum_{k=1}^{\sigma-1}\frac{a_k(\theta_x)}{t^{k}}\right)\Big|}{\frac{\exp(-t\eta_{\lambda,Z}^*(x))}{\sqrt{2\lambda\pi t\mathbb{E}[Z^2\mathrm{e}^{\theta_x Z}]}}}=+\infty,
\end{equation}
and let $\{t_k\}_{k\geq 1}\subset (0,\infty)$ be a divergent sequence that realizes the $\limsup$. By Lemma \ref{thm:modphiRspeed} the stochastic process
$\{S_{t_k}\}_{k\geq 1}$ converges mod-$\phi_{\lambda,Z}$ on $D_{cp}$ with speed $O(1/(t_k)^\sigma)$ and limiting function
\eqref{eq:psi}. So, by Theorem 3.2.2$(1)$ in \cite{FMN} easily follows that the relation \eqref{eq:deviationordertlattice} holds with $t_k$ in place of $t$.
Indeed, in our framework, the functions $F(\cdot)$, $h_{\cdot}$ and $\eta(\cdot)$ in the Theorem 3.2.2 of \cite{FMN} are given, respectively, by
the functions $\eta_{\lambda,Z}^*(\cdot)$ (in \eqref{eq:FenchelLegendre}), $\theta_{\cdot}$ and $\eta_{\lambda,Z}(\cdot)$. 
Moreover, the function $\psi(\cdot)$ and the quantities $c$ and $d$ of Theorem 3.2.2 
in \cite{FMN} are given, respectively, by the function 
$\psi(\cdot)$ defined by \eqref{eq:psi}, and $c=-\infty$, $d=a$. Therefore we reached a contradiction and the $\limsup$ in \eqref{eq:limsup7genSBR} is finite.
The proof is completed.\\
\noindent{\it Proof\,\,of\,\,Part\,\,$(ii)$.}\\
Similar to the proof of Part $(i)$. One has to apply Theorem 3.2.2$(2)$ in \cite{FMN} in place of Theorem 3.2.2$(1)$ in \cite{FMN}.  
\\
\noindent$\square$

\noindent{\it Proof\,\,of\,\,Lemma\,\,\ref{thm:modphiRspeed}.} Note that the assumptions imply \eqref{eq:Lq}. Indeed, for any $t>0$ and $q>1$,
\begin{equation}\label{eq:7gen}
\mathbb{E}\left[\left(\int_0^t(Z-H(s,M_1))\,\mathrm{d}s\right)^q\right]\leq t^q\mathbb{E}[Z^q]<\infty,
\end{equation}
where the latter quantity is finite by \eqref{eq:lighttailCP}. So \eqref{eq:Lq} follows by \eqref{eq:7gen}, \eqref{eq:Lqspeed} and Minkowski's inequality.
Consequently, all the claims in the proof of Lemma \ref{thm:modphiR} hold, and using the same notation of the lemma, we only have to prove that for any compact $K\subset D_{cp}$, there exists $C_K>0$ such that
\[
\sup_{z\in K}|\psi_t(z)-\psi(z)|\leq C_K t^{-\sigma}, 
\quad\text{for any $t>0$.}
\]
To this aim, we note that by the relations \eqref{eq:moduliexp1} and \eqref{eq:7gen1}, for any $t>0$ and $K\subset D_{cp}$, we have 
\begin{align}
|\psi_t(z)-\psi(z)|
&\leq\lambda\mathrm{e}^{\lambda\mathcal{M}_K(0)}\overline{M}_K\left(\mathbb{E}[\mathrm{e}^{p\,m_K Z}]^{1/p}\right)
\Big\|\int_t^\infty(Z-H(s,M_1))\,\mathrm{d}s\Big\|_{L^q(\mathbb P)}.\nonumber
\end{align}
The claim follows by the assumption \eqref{eq:Lqspeed} (note that the conjugate exponents $p$ and $q$ in 
this latter inequality depend on $K$, but not on $z$). 
\\
\noindent$\square$

\subsection{Computation of the quantities $a_{k}(\theta_x)$ and $b_k(\theta_x)$, $1\leq k\leq\sigma-1$, and proof of Proposition \ref{prop:estimateintegral}}
\label{sec:akprop33}

The following proposition provides the expression of the functions $a_k(\cdot)$ and $b_k(\cdot)$ in Theorem \ref{thm:deviationslattice}.

\begin{Proposition}\label{prop:akbk}
Let the assumptions and notation of Theorem \ref{thm:deviationslattice} prevail. Then,
for any $x\in (0,\lambda\mathbb{E}[Z\mathrm{e}^{a Z}])$ and $k=1,\ldots,\sigma-1$, we have
\begin{align}
a_k(\theta_x)&=\sum_{\ell=0}^{2k}\frac{\psi^{(2k-\ell)}(\theta_x)}{(2k-\ell)!}\sum_{(m_1,\ldots,m_\ell)\in\mathcal{S}_\ell}\frac{(-1)^{m_1+\ldots+m_\ell}}
{m_1!1!^{m_1}m_2!2!^{m_2}\ldots m_{\ell}!\ell!^{m_\ell}}\nonumber\\
&\qquad
\times
\prod_{j=1}^{\ell}
\left(\frac{\mathbb{E}[Z^{j+2}\mathrm{e}^{\theta_x Z}]}{(j+1)(j+2)\mathbb{E}[Z^2\mathrm{e}^{\theta_x Z}]}\right)^{m_j}
\frac{(-1)^{k}(2(k+m_1+\ldots+m_\ell)-1)!!}{(\lambda\mathbb{E}[Z^2\mathrm{e}^{\theta_x Z}])^k}\nonumber
\end{align}
and
\begin{align}
b_k(\theta_x)&=\sum_{n=0}^{2k}\sum_{(m_1,\ldots,m_n)\in\mathcal{S}_n}\frac{\mathrm{e}^{(m_1+\ldots+m_n)\theta_x}
(m_1+\ldots+m_n)!(1-\mathrm{e}^{-\theta_x})^{-(m_1+\ldots+m_n)-1}}{m_1!1!^{m_1}m_2!2!^{m_2}\ldots m_{n}!n!^{m_n}}
\prod_{j=1}^{n}(-1)^{jm_j}\nonumber\\
&\qquad
\times\sum_{\ell=0}^{2k-n}\frac{\psi^{(2k-n-\ell)}(\theta_x)}{(2k-n-\ell)!}\sum_{\mathcal{S}_\ell}\frac{(-1)^{m_1+\ldots+m_\ell}}
{m_1!1!^{m_1}m_2!2!^{m_2}\ldots m_{\ell}!\ell!^{m_\ell}}\nonumber\\
&\qquad
\times\prod_{j=1}^{\ell}
\left(\frac{\mathbb{E}[Z^{j+2}\mathrm{e}^{\theta_x Z}]}{(j+1)(j+2)\mathbb{E}[Z^2\mathrm{e}^{\theta_x Z}]}
\right)^{m_j}
\frac{(-1)^{k}(2(k+m_1+\ldots+m_\ell)-1)!!}{(\lambda\mathbb{E}[Z^2\mathrm{e}^{\theta_x Z}])^k}.\nonumber
\end{align}
Here $\mathcal{S}_n$ denotes the set of $n$tuples of non-negative integers $(m_1,\ldots,m_n)$ such that
$1\cdot m_1+2\cdot m_2+\ldots n\cdot m_n=n$, $f^{(k)}$ denotes the derivative of order $k\in\mathbb N$ of a sufficiently smooth 
function $f$ and $f^{(0)}=f$. The derivatives of the function $\psi(\cdot)$, defined by \eqref{eq:psi}, can be computed by using the recursive formula:
\[
\psi^{(n)}(\theta_x)=n!\sum_{i=1}^{n}\lambda^i\frac{\mathrm{e}^{\lambda\varphi(\theta_x)}}{i!}
\sum_{m_1+m_2+\ldots+m_i=n}\frac{\varphi^{(m_1)}(\theta_x)}{m_1!}\ldots\frac{\varphi^{(m_i)}(\theta_x)}{m_i!},
\quad n\in\mathbb N,
\]
where the sum is taken over all the $m_1,\ldots,m_i\in\mathbb N$ such that $m_1+\ldots+m_i=n$, and
\[
\varphi^{(n)}(\theta_x)=\int_0^\infty(\mathbb{E}[H(s,M_1)^n \mathrm{e}^{\theta_xH(s,M_1)}]-\mathbb{E}[Z^n\mathrm{e}^{\theta_x Z}])\,\mathrm{d}s.
\]
\end{Proposition}

\noindent{\it Proof.} The general expression of the functions $a_k(\cdot)$ and $b_k(\cdot)$ can be computed following the suggestions of Remark 3.2.5
in \cite{FMN}. See Proposition 1 in \cite{GZ} for the details. The recursive formula for the derivatives of $\psi$ is easily 
obtained by applying the Fa\`a di Bruno formula (see Lemma \ref{le:FaadiBruno}). \\
\noindent$\square$

\noindent{\it Proof\,\,of\,\,Proposition\,\,\ref{prop:estimateintegral}.}\\
\noindent{\it Proof\,\,of\,\,Part\,\,$(i)$.}\\
The lower bound is trivial. The upper bound follows immediately by applying Lemma \ref{le:standardineq}.\\
\noindent{\it Proof\,\,of\,\,Part\,\,$(ii)$.}\\
The upper bound is trivial. For the lower bound, note that by Lemma \ref{le:standardineq} and H\"older's inequality it follows
\begin{align}
\int_0^\infty(\mathbb{E}[\mathrm{e}^{\theta_x Z}]-\mathbb{E}[\mathrm{e}^{\theta_x H(s,M_1)}])\,\mathrm{d}s
&\leq\theta_x\|\mathrm{e}^{\theta_x Z}\|_{L^p(\mathbb P)}
\Big\|\int_0^\infty(Z-H(s,M_1))\,\mathrm{d}s\Big\|_{L^q(\mathbb P)}.\nonumber
\end{align}
\noindent$\square$

\subsection{Proof of Theorem \ref{thm:fluctuations}}\label{sec:flucFeb2}

\noindent{\it Proof\,\,of\,\,Part\,\,$(i)$.}\\ 
Reasoning by contradiction, suppose
\begin{equation}\label{eq:2Feb2021}
\limsup_{t\to\infty}\frac{\mathbb{P}\left(\frac{S_t-\lambda\mathbb{E}[Z]t}{\sqrt{\lambda\mathbb{E}[Z^2]t}}>y(t)\right)}
{\mathbb{P}(\mathrm{N}(0,1)>y(t))}\neq 1,
\end{equation}
and let $\{t_k\}_{k\geq 1}\subset (0,\infty)$ be a divergent sequence that realizes the $\limsup$. The proof of \eqref{eq:CLT} then proceeds similarly to the proof of 
\eqref{eq:deviationordert} in Theorem \ref{thm:deviationsnonlattice} [i.e., applying Lemma \ref{thm:modphiR} and then
Theorem 4.3.1 in \cite{FMN} (in place of Theorem 4.2.1 in \cite{FMN}) to reach a contradiction]. By a similar argument we reach a contradiction if we replace
the $\limsup$ with the $\liminf$ in \eqref{eq:2Feb2021}.
\\ 
\noindent{\it Proof\,\,of\,\,Part\,\,$(ii)$.}\\ 
{From now on, for ease of notation,  we denote by $f$ the Fenchel-Legendre transform $\eta_{\lambda,Z}^{*}$ in 
\eqref{eq:FenchelLegendre}.  For technical reasons (related to the application of Dini's implicit function theorem), we assume  $f$ to be
defined in a neighborhood of $\lambda\mathbb{E}[Z]$}. The relations in \eqref{eq:extCLT} follow
(reasoning again by contradiction as in the Part $(i)$) by Lemma \ref{thm:modphiR} and 
Theorem 4.3.1 in \cite{FMN}. As far as \eqref{eq:Dic3mattino1} is concerned, 
note that by the definition of the function $\theta_{\cdot}$ and Dini's implicit function theorem we have that $\theta_{\cdot}$ is infinitely differentiable on 
the domain of $f$ and
\[
\theta'_x=\frac{1}{\lambda\mathbb{E}[Z^2\mathrm{e}^{\theta_x Z}]}.
\]
In particular, $f$ is infinitely differentiable on its domain and,
by Taylor's formula, in a neighborhood of $\lambda\mathbb{E}[Z]$, we have
\[
f(x)=f(\lambda\mathbb{E}[Z])+f'(\lambda\mathbb{E}[Z])(x-\lambda\mathbb{E}[Z])+\frac{f''(\lambda\mathbb{E}[Z])}{2}(x-\lambda\mathbb{E}[Z])^2
+o((x-\lambda\mathbb{E}[Z])^2).
\]
Note that $\theta_{\lambda\mathbb{E}[Z]}=0$ and so $f(\lambda\mathbb{E}[Z])=0$.
A straightforward computation yields
\[
f'(x)=\theta_{x}+x\theta'_x-\lambda\theta'_x\mathbb{E}[Z\mathrm{e}^{\theta_{x}Z}]=\theta_x\quad\text{and}\quad f''(x)=\theta'_x.
\]
Therefore,
\[
f'(\lambda\mathbb{E}[Z])=\theta_{\lambda\mathbb{E}[Z]}=0\quad\text{and}\quad
f''(\lambda\mathbb{E}[Z])=\theta'_{\lambda\mathbb{E}[Z]}=\frac{1}{\lambda\mathbb{E}[Z^2]}.
\]
So, as $x\to\lambda\mathbb{E}[Z]$, we have
\[
f(x)=\frac{1}{2\lambda\mathbb{E}[Z^2]}(x-\lambda\mathbb{E}[Z])^2+
o((x-\lambda\mathbb{E}[Z])^2).
\]
In particular, as $t\to\infty$, by the expression of $v(\cdot)$ we have
\[
f(v(t))=\frac{y(t)^2}{2t}+o\left(\frac{y(t)^2}{t}\right),
\]
which gives the claim.
\\
\noindent{\it Proof\,\,of\,\,Part\,\,$(iii)$.}\\ 
Note that, in particular, $y(t)=o(t^{1/2})$, and therefore by the previous Part $(ii)$ we have that the
relations in \eqref{eq:extCLT} hold. 
In the previous Part $(ii)$ we noticed that $f$ is infinitely differentiable
on its domain 
with $f''(\cdot)=\theta'_{\cdot}$. Therefore, the $k$th derivative of $f$ is given by $f^{(k)}(x)=\theta_x^{(k-1)}$, $k\geq 1$. By Taylor's formula, as $x\to\lambda\mathbb{E}[Z]$, we have
\[
f(x)=\sum_{k=2}^{m}\frac{\theta_{\lambda\mathbb{E}[Z]}^{(k-1)}}{k!}(x-\lambda\mathbb{E}[Z])^{k}+o((x-\lambda\mathbb{E}[Z])^m).
\]
In particular, as $t\to\infty$, by the expression of $v(\cdot)$ we have
\begin{equation}\label{eq:Dic3seraIII}
tf(v(t))=\sum_{k=2}^{m}(\lambda\mathbb{E}[Z^2])^{k/2}
\frac{\theta_{\lambda\mathbb{E}[Z]}^{(k-1)}}{k!}\frac{y(t)^k}{t^{k/2-1}}+t\cdot o\left(\frac{y(t)^m}{t^{m/2}}\right).
\end{equation}
Note that 
\begin{align}
\lim_{t\to\infty}t\cdot o\left(\frac{y(t)^m}{t^{m/2}}\right)
&=\lim_{t\to\infty}\frac{o\left(\frac{y(t)^m}{t^{m/2}}\right)}{\frac{y(t)^m}{t^{m/2}}}\lim_{t\to\infty}\left(\frac{y(t)}{t^{\frac{1}{2}-\frac{1}{m}}}\right)^m=0,\label{eq:Dic3seraII}
\end{align}
where we used the assumption on $y(\cdot)$. Moreover,
\begin{align}
\sum_{k=2}^{m}(\lambda\mathbb{E}[Z^2])^{k/2}
\frac{\theta_{\lambda\mathbb{E}[Z]}^{(k-1)}}{k!}\frac{y(t)^k}{t^{k/2-1}}
&=y(t)^2\left(\frac{1}{2}+\sum_{k=3}^{m}(\lambda\mathbb{E}[Z^2])^{k/2}
\frac{\theta_{\lambda\mathbb{E}[Z]}^{(k-1)}}{k!}\frac{y(t)^{k-2}}{t^{k/2-1}}\right)\nonumber\\
&=y(t)^2\left(\frac{1}{2}+\sum_{j=1}^{m-2}(\lambda\mathbb{E}[Z^2])^{(j+2)/2}
\frac{\theta_{\lambda\mathbb{E}[Z]}^{(j+1)}}{(j+2)!}\left(\frac{y(t)}{\sqrt t}\right)^j\right)\to\infty,\quad\text{as $t\to\infty$.}\label{eq:Dic3seraI}
\end{align}
The claim follows combining \eqref{eq:extCLT} with
\eqref{eq:Dic3seraIII}, \eqref{eq:Dic3seraII} and \eqref{eq:Dic3seraI}.
Finally, the recursive formula for $\theta_{\lambda\mathbb{E}[Z]}^{(j)}$,
$2\leq j\leq m-1$, easily follows by the Fa\`a di Bruno formula (see Lemma \ref{le:FaadiBruno}). The proof is completed.

\subsection{Proof of Propositions \ref{prop:PoissonCluster} and \ref{prop:GW}}

\noindent{\it Proof\,\,of\,\,Proposition\,\,\ref{prop:PoissonCluster}.}\\
\noindent{\it Proof\,\,of\,\,Parts\,\,$(i)$\,\,and\,\,$(ii)$.}\\
The claim follows by applying Theorems \ref{thm:deviationslattice} and \ref{thm:fluctuations}. To this aim, we show that, under the foregoing assumptions,
one has \eqref{eq:Lqspeed} (then, due to the fact that $Z:=N_1([0,\infty))$ satisfies \eqref{eq:lighttailCP},  one also has that condition
\eqref{eq:Lq} holds). We start noticing that, for $t\geq 0$, we have
\begin{align}
\int_{t}^{\infty}(N_1([0,\infty))-N_1([0,s]))\,\mathrm{d}s&=\int_{t}^{\infty}(N_1([0,L])-N_1([0,s]))\,\mathrm{d}s\nonumber\\
&=\bold{1}(L>t)\int_{t}^{\infty}(N_1([0,L])-N_1([0,s]))\bold{1}(L>s)\,\mathrm{d}s\nonumber\\
&\leq (L-t)^{+}N_1((t,L]),\label{eq:ineqvaccino}
\end{align}
where $x^+:=\max\{x,0\}$.  For an arbitrary $q>1$,  set $\widetilde{p}_1:=\frac{m}{q}$, $m>q$, $m\in\mathbb N$, and let $\widetilde{p}_2:=m/(m-q)$. Note that $\widetilde{p}_1,\widetilde{p}_2>1$ and 
$(\widetilde{p}_1)^{-1}+(\widetilde{p}_2)^{-1}=1$. Therefore by H\"older's inequality we have
\begin{align}
\Big\|\int_{t}^{\infty}(N_1([0,L])-N_1([0,s]))\,\mathrm{d}s\Big\|_{L^q(\mathbb P)}
&\leq\mathbb{E}[((L-t)^{+})^{m}]^{1/m}\mathbb{E}[N_1((t,L])^{q\widetilde{p}_2}]^{1/(q\widetilde{p}_2)}.\nonumber
\end{align}
For $\theta^*\in (0,a)$ fixed, we have 
\[
\mathbb{E}[N_1((t,L])^{q\widetilde{p}_2}]^{1/(q\widetilde{p}_2)}\leq\mathbb{E}[(N_1([0,\infty))^{\lceil q\widetilde{p}_2\rceil}]\leq\frac{(\lceil q\widetilde{p}_2\rceil)!}{(\theta^*)^{\lceil q\widetilde{p}_2\rceil}}
\mathbb{E}[\mathrm{e}^{\theta^{*}N_1([0,\infty))}]=:\kappa_q^{(1)},
\]
where $x\mapsto\lceil x\rceil$ denotes the ceiling function. So
\begin{align}
\Big\|\int_{t}^{\infty}(N_1([0,L])-N_1([0,s]))\,\mathrm{d}s\Big\|_{L^q(\mathbb P)}
&\leq\kappa_q^{(1)}
\mathbb{E}[((L-t)^{+})^{m}]^{1/m}.\label{eq:AI}
\end{align}
Let $\sigma\in\mathbb N$ be arbitrarily fixed and set $\theta:=\sigma m+1$.
Since $(L-t)^+$ is a non-negative random variable,  for any $t\geq 1$ by Markov's inequality we have
\begin{align}
\mathbb{E}[((L-t)^{+})^{m}]&=m\int_0^\infty u^{m-1}\mathbb{P}((L-t)^+>u)\,\mathrm{d}u\nonumber\\
&=m\int_0^\infty u^{m-1}\mathbb{P}(L^{\theta m}>(u+t)^{\theta m})\,\mathrm{d}u\nonumber\\
&\leq m\mathbb{E}[L^{\theta m}]\int_0^\infty\frac{u^{m-1}}{(u+t)^{\theta m}}\,\mathrm{d}u\nonumber\\
&=\mathbb{E}[L^{\theta m}]\int_0^\infty\frac{1}{(z^{1/m}+t)^{\theta m}}\,\mathrm{d}z\nonumber\\
&\leq\mathbb{E}[L^{\theta m}]\int_0^\infty\frac{1}{z^{\theta}+t^{\theta}}\,\mathrm{d}z\nonumber\\
&=\mathbb{E}[L^{\theta m}]t^{-(\theta-1)}\int_0^\infty\frac{1}{z^{\theta}+1}\,\mathrm{d}z\nonumber\\
&\leq\mathbb{E}[L^{\theta m}]t^{-(\theta-1)}\left(1
+\int_1^\infty\frac{1}{z^{2}+1}\,\mathrm{d}z\right)
:=\kappa_q^{(2)}t^{-(\theta-1)},\nonumber
\end{align}
where the quantity which multiplies $t^{-(\theta-1)}$ depends on $q$ since $m=m(q)$.  So
\[
\mathbb{E}[((L-t)^{+})^{m}]^{1/m}\leq(\kappa_q^{(2)})^{1/m}t^{-(\theta-1)/m}\leq (\kappa_q^{(2)})^{1/m}\,t^{-\sigma},
\]
where we used the definition of $\theta$ and the fact that $t\geq 1$.
Combining this latter inequality with \eqref{eq:AI}, for any $\sigma\in\mathbb N$ fixed and $t\geq 1$, we have
\[
\Big\|\int_{t}^{\infty}(N_1([0,L])-N_1([0,s]))\,\mathrm{d}s\Big\|_{L^q(\mathbb P)}
\leq\kappa_q t^{-\sigma},
\]
where $\kappa_q:=\kappa_q^{(1)}(\kappa_q^{(2)})^{1/m}$. 
Therefore, for any $\sigma\in\mathbb{N}$ fixed we have
\[
\sup_{q>1}\kappa_q^{-1}\Big\|\int_{t}^{\infty}(N_1([0,L])-N_1([0,s]))\,\mathrm{d}s\Big\|_{L^q(\mathbb P)}
\leq t^{-\sigma},\quad\text{for all $t\geq 1$.}
\]
Thus condition \eqref{eq:Lqspeed} holds, and the proof of $(i)$ and $(ii)$ is completed. \\
\noindent{\it Proof\,\,of\,\,Part\,\,$(iii)$.}\\
The claim easily follows by Proposition \ref{prop:estimateintegral}. Here, we limit ourselves to note that
by \eqref{eq:ineqvaccino} we have
\[
\int_{0}^{\infty}(Z-N_1([0,s]))\,\mathrm{d}s\leq LZ,
\]
and so the upper bound in \eqref{eq:lowerbd2800} and the lower bound in \eqref{eq:lowerbd2801} follow, respectively, 
by the upper bound in the Part $(i)$ of Proposition \ref{prop:estimateintegral} and the lower bound in the Part $(ii)$ of Proposition \ref{prop:estimateintegral},
and H\"older's inequality.\\
\noindent$\square$
\\
\noindent{\it Proof\,\,of\,\,Proposition\,\,\ref{prop:GW}.}\\
\noindent{\it Proof\,\,of\,\,Parts\,\,$(i)$\,\,and\,\,$(ii)$.}\\
We are going to apply Proposition \ref{prop:PoissonCluster}, and so we verify the assumptions therein.
We already noticed that $Z:=N_1([0,\infty))$ satisfies \eqref{eq:lighttailCP}.
So we only need to check \eqref{eq:Lpiccolo}. 
Let $\{U_n\}_{n\geq 1}$ be a sequence of independent random variables with the same law as $B_{1,1,1}$, independent of $Z$.
Since $V$, defined by \eqref{eq:V}, has the same law as
\[
\sum_{n=1}^{Z-1}U_n,
\]
and $V$ is an upper bound of $L$,  we have
\begin{equation}\label{eq:IUB2901}
L\leq L':=\sum_{n=1}^{Z}U_n.
\end{equation}
Such random variable $L'$ has all the moments finite, indeed
for any $k\in\mathbb N$, using Minkowski's inequality we have
\begin{equation}\label{eq:IIUB2901}
\mathbb{E}\left[\left(\sum_{n=1}^{Z}U_n\right)^k\right]\leq\mathbb{E}[Z^k]\mathbb{E}[U_1^k]<\infty,
\end{equation}
where the latter term is finite due to \eqref{hp:densityCramer} and the fact that $Z$
satisfies \eqref{eq:lighttailCP}. The proof is completed.\\
\noindent{\it Proof\,\,of\,\,Part\,\,$(iii)$.}\\
By the random walk perspective to branching processes (see e.g. \cite{Remco} p. 90), we have that $Z$ has the same law of the
hitting time to zero of the random walk $\{Z'_n\}_{n\geq 0}$ defined by
\[
Z'_0:=1\quad Z'_n:=P_{1,1}+\ldots+P_{1,n}-(n-1),\quad n\in\mathbb N,
\]
i.e.
\[
Z\overset{d}{=}\inf\{n\in\mathbb N:\,\,Z'_n=0\}
\]
(we define the right-hand side equal to $+\infty$ if $\{\ldots\}=\emptyset$). So
\begin{align}
\mathbb{P}(Z>n)\leq\mathbb{P}(Z'_n>0)
&=\mathbb{P}\left(\frac{1}{n}\sum_{k=1}^{n}P_{1,k}\geq 1\right)
\leq\mathrm{e}^{-n a_c},\quad n\in\mathbb N,\label{eq:tailtotprog}
\end{align}
where the latter inequality follows by the Chernoff's bound (see e.g. \cite{Remco} Theorem 2.19 p. 68, and recall that $\mathbb{E}[P_{1,1}]<1$).
Standard convexity arguments (combined with the subcritical assumption) guarantee $a_c\in (0,\infty)$.  A straightforward computation shows that
\eqref{eq:tailtotprog} implies $\mathbb{E}[\mathrm{e}^{\theta Z}]<\infty$ $\forall$ $\theta<a_c$, which completes the proof.\\
\noindent{\it Proof\,\,of\,\,Part\,\,$(iv)$.}\\
Let $\{Z_j\}_{j\geq 1}$ be independent copies of $Z$ and let $K_1$ be the random variable defined at the beginning of Subsection \ref{subsec:Kuno}.
For any $\theta\in(-\infty,b_c)$, 
by standard computations we have
\begin{align}
\mathbb{E}[\mathrm{e}^{\theta Z}]
&=\mathrm{e}^{\theta}\sum_{k\geq 0}\mathbb{E}[\mathrm{e}^{\theta\sum_{j=1}^{k}Z_j}\,|\,K_1=k]p_k\nonumber\\
&=\mathrm{e}^{\theta}\sum_{k\geq 0}\mathbb{E}[\mathrm{e}^{\theta Z}]^{k}p_k=\mathrm{e}^{\theta}G_{P_{1,1}}(\mathbb{E}[\mathrm{e}^{\theta Z}])<\infty.
\label{eq:eqLap}
\end{align}
As far as the moments of $Z$ are concerned, 
we combine the Fa\`a di Bruno formula with the elementary relation: 
\begin{equation}\label{eq:fprod}
\frac{\mathrm{d}^n}{\mathrm{d}x^n}(f(x)g(x))=\sum_{k=0}^{n}\binom{n}{k}f^{(n-k)}(x)g^{(k)}(x),\quad n\in\mathbb N,
\end{equation}
for sufficiently smooth functions $f$ and $g$. For any $\theta\in (-\infty,b_c)$, by \eqref{eq:eqLap} and \eqref{eq:fprod} we have
\[
\frac{\mathrm{d}^n}{\mathrm{d}\theta^n}\mathbb{E}[\mathrm{e}^{\theta Z}]=\mathrm{e}^{\theta}\sum_{k=0}^{n}\binom{n}{k}
\frac{\mathrm{d}^k}{\mathrm{d}\theta^k}G_{P_{1,1}}(\mathbb{E}[\mathrm{e}^{\theta Z}]).
\]
By the Fa\`a di Bruno formula we have
\[
\frac{\mathrm{d}^k}{\mathrm{d}\theta^k}G_{P_{1,1}}(\mathbb{E}[\mathrm{e}^{\theta Z}])
=k!\sum_{i=1}^{k}\frac{G_{P_{1,1}}^{(i)}(\mathbb{E}[\mathrm{e}^{\theta Z}])}{i!}
\sum_{m_1+m_2+\ldots+m_i=k}\frac{
\frac{\mathrm{d}^{m_1}}{\mathrm{d}\theta^{m_1}}
\mathbb{E}[\mathrm{e}^{\theta Z}]}{m_1!}\ldots\frac{\frac{\mathrm{d}^{m_i}}{\mathrm{d}\theta^{m_i}}
\mathbb{E}[\mathrm{e}^{\theta Z}]}{m_i!},
\quad k\in\mathbb N,
\]
where the sum is taken over all the $m_1,\ldots,m_i\in\mathbb N$ such that $m_1+\ldots+m_i=k$. 
The claim follows recalling that
\[
\mathbb{E}[Z^n]=\frac{\mathrm{d}^n}{\mathrm{d}\theta^n}\mathbb{E}[\mathrm{e}^{\theta Z}]\Big|_{\theta=0}
\]
and that
\[
G_{P_{1,1}}^{(i)}(1)=\mathbb{E}[P_{1,1}(P_{1,1}-1)\ldots(P_{1,1}-(i-1))]
\]
is  the $i$th factorial moment of $P_{1,1}$.\\
\noindent{\it Proof\,\,of\,\,Part\,\,$(v)$.}\\
Differentiating \eqref{eq:eqLap} with respect to $\theta$, for any $\theta\in (-\infty,b_c)$ we have
\begin{align}
\mathbb{E}[Z\mathrm{e}^{\theta Z}]
&=\mathrm{e}^{\theta}G_{P_{1,1}}(\mathbb{E}[\mathrm{e}^{\theta Z}])+\mathrm{e}^{\theta}\mathbb{E}[Z\mathrm{e}^{\theta Z}]
G'_{P_{1,1}}(\mathbb{E}[\mathrm{e}^{\theta Z}])
=\mathbb{E}[\mathrm{e}^{\theta Z}]+\mathrm{e}^{\theta}\mathbb{E}[Z\mathrm{e}^{\theta Z}]
G'_{P_{1,1}}(\mathbb{E}[\mathrm{e}^{\theta Z}]),\nonumber
\end{align}
where the latter equality follows by \eqref{eq:eqLap}.  Therefore
\begin{align}
\mathbb{E}[Z\mathrm{e}^{\theta Z}]
&=\frac{\mathbb{E}[\mathrm{e}^{\theta Z}]}{1-\mathrm{e}^{\theta}G'_{P_{1,1}}(\mathbb{E}[\mathrm{e}^{\theta Z}])}.\label{eq:22Gen}
\end{align}
Differentiating again we have
\begin{equation}\label{eq:23Gen}
\mathbb{E}[Z^2\mathrm{e}^{\theta Z}]=\mathbb{E}[\mathrm{e}^{\theta Z}]
\frac{1+\mathrm{e}^{\theta}(G'_{P_{1,1}}(\mathbb{E}[\mathrm{e}^{\theta Z}])+
G''_{P_{1,1}}(\mathbb{E}[\mathrm{e}^{\theta Z}]))}
{(1-\mathrm{e}^{\theta}G'_{P_{1,1}}(\mathbb{E}[\mathrm{e}^{\theta Z}]))^2}.
\end{equation}
Recall that, for any $x\in (0,\lambda\mathbb{E}[Z\mathrm{e}^{b_c Z}])$ we denote by $\theta_x\in (-\infty,b_c)$ the unique solution to
$\lambda\mathbb{E}[Z\mathrm{e}^{\theta Z}]=x$.  Using \eqref{eq:eqLap}, we rewrite the relation \eqref{eq:22Gen} as
\begin{align*}
\mathbb{E}[Z\mathrm{e}^{\theta Z}]
&=\frac{\mathbb{E}[\mathrm{e}^{\theta Z}]}{1-\mathbb{E}[\mathrm{e}^{\theta Z}]\frac{G'_{P_{1,1}}(\mathbb{E}[\mathrm{e}^{\theta Z}])}{G_{P_{1,1}}(\mathbb{E}[\mathrm{e}^{\theta Z}])}}.\nonumber
\end{align*}
Taking $\theta=\theta_x$ in this relation, we have
\[
\frac{x}{\lambda}
=\frac{\mathbb{E}[\mathrm{e}^{\theta_x Z}]}{1-\mathbb{E}[\mathrm{e}^{\theta_x Z}]\frac{G'_{P_{1,1}}(\mathbb{E}[\mathrm{e}^{\theta_x Z}])}{G_{P_{1,1}}(\mathbb{E}[\mathrm{e}^{\theta_x Z}])}},\nonumber
\]
and so $\varrho_x=\mathbb{E}[\mathrm{e}^{\theta_x Z}]$ is solution to \eqref{eq:equationSab}. The expressions of $\theta_x$
and $\mathbb{E}[Z^2\mathrm{e}^{\theta_x Z}]$ are readily obtained by \eqref{eq:eqLap} and \eqref{eq:23Gen}, respectively.
The proof is completed.\\
\noindent{\it Proof\,\,of\,\,Part\,\,$(vi)$.}\\
For a fixed $r>0$, let $\{N_{r}^{(j)}(\cdot)\}_{j\geq 1}$ be independent and identically distributed point processes  on $[r,\infty)$
with the same branching structure as $N_1$, but ancestor in $r$
(when $r=0$ they are independent copies of $N_1$). Hereon, we denote by 
$N_r^{(j)}(s)$ the number of points of $N_r^{(j)}$ on the interval $[r,s]$ and by $Q$ the law of $B_{1,1,1}$.
For any $\theta\in(-\infty,b_c)$ and $s>0$,  similarly to
\eqref{eq:eqLap}, we have
\begin{align}
\mathbb{E}[\mathrm{e}^{\theta N_1([0,s])}]&
=\sum_{k\geq 0}\mathbb{E}[\mathrm{e}^{\theta N_1([0,s])}\,|\,K_1=k]p_k\nonumber\\
&=\sum_{k\geq 0}\mathbb{E}\left[\mathrm{e}^{\theta\left(1+\sum_{j=1}^{k}N_{B_{1,0,j}}^{(j)}(s)\right)}\,|\,K_1=k\right]p_k\nonumber\\
&=\mathrm{e}^{\theta}\sum_{k\geq 0}p_k\left(\int_{[0,\infty)}\mathbb{E}[\mathrm{e}^{\theta N_r^{(1)}(s)}]Q(\mathrm{d}r)\right)^k\nonumber\\
&=\mathrm{e}^{\theta}G_{P_{1,1}}\left(\int_{[0,\infty)}\mathbb{E}[\mathrm{e}^{\theta N_r^{(1)}(s)}]Q(\mathrm{dr})\right).\nonumber
\end{align}
By this relation, \eqref{eq:eqLap} and the mean value theorem,  we have (since $N_r^{(1)}(s)=0$ for $r>s$)
\begin{align}
\mathbb{E}[\mathrm{e}^{\theta N_1([0,s])}]-\mathbb{E}[\mathrm{e}^{\theta Z}]
&=\mathrm{e}^{\theta}\left(G_{P_{1,1}}\left(\int_{[0,\infty)}\mathbb{E}[\mathrm{e}^{\theta N_r^{(1)}(s)}]Q(\mathrm{dr})\right)-G_{P_{1,1}}(\mathbb{E}[\mathrm{e}^{\theta Z}])\right)\nonumber\\
&=\mathrm{e}^{\theta}\left(\int_{[0,\infty)}\mathbb{E}[\mathrm{e}^{\theta N_r^{(1)}(s)}]Q(\mathrm{dr})-\mathbb{E}[\mathrm{e}^{\theta Z}]\right)G'_{P_{1,1}}(\xi)\nonumber\\
&=\mathrm{e}^{\theta}\left(\int_{[0,s]}(\mathbb{E}[\mathrm{e}^{\theta N_r^{(1)}(s)}]-\mathbb{E}[\mathrm{e}^{\theta Z}])Q(\mathrm{dr})+(1-\mathbb{E}[\mathrm{e}^{\theta Z}])
Q((s,\infty))\right)G'_{P_{1,1}}(\xi),\nonumber
\end{align}
for some
\[
\xi\in\left(\min\left\{\int_{[0,\infty)}\mathbb{E}[\mathrm{e}^{\theta N_r^{(1)}(s)}]Q(\mathrm{dr}),\mathbb{E}[\mathrm{e}^{\theta Z}]\right\},
\max\left\{\int_{[0,\infty)}\mathbb{E}[\mathrm{e}^{\theta N_r^{(1)}(s)}]Q(\mathrm{dr}),\mathbb{E}[\mathrm{e}^{\theta Z}]\right\}
\right).
\]
Since $\xi\geq 0$, a simple computation shows that $G'_{P_{1,1}}(\xi)\geq p_1$. Therefore
\begin{align}
\mathbb{E}[\mathrm{e}^{\theta N_1([0,s])}]-\mathbb{E}[\mathrm{e}^{\theta Z}]
&\leq\mathrm{e}^{\theta}\left(\int_{[0,s]}(\mathbb{E}[\mathrm{e}^{\theta N_r^{(1)}(s)}]-\mathbb{E}[\mathrm{e}^{\theta Z}])Q(\mathrm{dr})
+(1-\mathbb{E}[\mathrm{e}^{\theta Z}])Q((s,\infty))\right)p_1\nonumber\\
&\leq p_1\mathrm{e}^{\theta}(1-\mathbb{E}[\mathrm{e}^{\theta Z}])
Q((s,\infty)),\quad\text{for any $\theta\in (0, b_c)$,}\nonumber
\end{align}
and
\begin{align}
\mathbb{E}[\mathrm{e}^{\theta N_1([0,s])}]-\mathbb{E}[\mathrm{e}^{\theta Z}]
\geq p_1\mathrm{e}^{\theta}(1-\mathbb{E}[\mathrm{e}^{\theta Z}])
Q((s,\infty)),\quad\text{for any $\theta<0$.}\nonumber
\end{align}
The claim follows combining these inequalities with \eqref{eq:lowerbd2800}, \eqref{eq:lowerbd2801},
\eqref{eq:IUB2901} and \eqref{eq:IIUB2901}.\\
\noindent$\square$

\subsection{Proof of Proposition \ref{cor:improvement}}

We only prove Part $(i)$. Indeed, mutatis mutandis (i.e.  applying Theorem \ref{thm:deviationslattice}$(ii)$ in place of
Theorem \ref{thm:deviationsnonlattice}) the proof of Part $(ii)$ is 
similar. We start noticing that the upper bound is a simple consequence of the inequality $\psi_{IBNR}(u)\leq\psi_{CL}(u)$, $u>0$, and
\eqref{eq:CLestimate}. As far as the lower bound is concerned, note that, for any $u,d>0$, we have
\begin{equation}\label{eq:SAmbrogio1}
\psi_{IBNR}(u)\geq\mathbb{P}(S_{ud}-cud\geq u)=\mathbb{P}\left(\frac{S_{ud}}{ud}\geq c+\frac{1}{d}\right).
\end{equation}
By the strict convexity of the function
\[
(0,a)\ni\gamma\mapsto\lambda(\mathbb{E}[\mathrm{e}^{\gamma Z}]-1)-c\gamma,
\]
we have $\lambda\mathbb{E}[Z\mathrm{e}^{w Z}]-c>0$. 
Set $d:=(\lambda\mathbb{E}[Z\mathrm{e}^{w Z}]-c)^{-1}>0$ 
and $x:=c+d^{-1}=\lambda\mathbb{E}[Z\mathrm{e}^{w Z}]$. By the \lq\lq net profit" condition and the fact that
$\mathbb{E}[Z\mathrm{e}^{a Z}]>\mathbb{E}[Z\mathrm{e}^{wZ}]>0$, we have $x\in (\lambda\mathbb{E}[Z],\lambda\mathbb{E}[Z\mathrm{e}^{aZ}])$. 
Moreover, a straightforward computation 
gives
\[
\theta_x=w\quad\text{and}\quad x\theta_x-\lambda(\mathbb{E}[\mathrm{e}^{\theta_x Z}]-1)=w(\lambda\mathbb{E}[Z\mathrm{e}^{wZ}]-c).
\]
The lower bound follows by these relations, \eqref{eq:SAmbrogio1} and Theorem \ref{thm:deviationsnonlattice}.

\subsection{Proof of Theorems \ref{cor:coninlaw} and \ref{prop:speed}}\label{sec:stable}

In this section we prove Theorems \ref{cor:coninlaw} and \ref{prop:speed}.  In particular, we emphasize that
the proof of Theorem \ref{prop:speed}
exploits the ideas and the techniques developed in \cite{FMN2}. 

\noindent{\it Proof of Theorem \ref{cor:coninlaw}.} We divide the proof in two steps. In the first step we prove the inequality
\begin{equation}\label{eq:2021Feb15}
|\psi_t(\bold{i}\xi)-1|\leq(H_1(t)|\xi|^{\alpha}+H_2(t)|\xi|^{2\alpha})
\exp(H_1(t)|\xi|^{\alpha}+H_2(t)|\xi|^{2\alpha}),\quad\text{$\forall$ $(t,\xi)\in (0,\infty)\times\mathbb{R}$,}
\end{equation}
where
\begin{equation}\label{eq:tetat}
\psi_t(\bold{i}\xi):=\mathbb{E}[\mathrm{e}^{\bold{i}\xi(S_t/t^{1/(2\alpha)})}]\mathrm{e}^{-\sqrt{t}\eta_{c\lambda^{1/\alpha},\alpha,\beta}(\bold{i}\xi)},\quad\text{$t>0$, $\xi\in\mathbb{R}$,}
\end{equation}
\begin{equation}\label{eq:H1primo}
H_1(t):=\lambda c^{\alpha}\sqrt{1+(\beta\tan(\pi\alpha/2))^2}\,\frac{\int_0^t(1-F(s)^{\alpha})\,\mathrm{d}s}{\sqrt t},\quad t>0,
\end{equation}
and
\begin{equation}\label{eq:K2primo}
H_2(t):=\frac{\lambda c^{2\alpha}}{2}[1+(\beta\tan(\pi\alpha/2))^2],\quad t>0.
\end{equation}
In the second step we conclude the proof.
\\
\noindent{\it Step\,\,1:\,\,Proof\,\,of\,\,\eqref{eq:2021Feb15}.}\\
For any $t>0$ and $\xi\in\mathbb R$,  by Lemma \ref{le:standardineq}
and the expression of the Laplace functional of a Poisson process (see e.g. \cite{DVJ}),
we have
\begin{align}
\mathbb{E}[\mathrm{e}^{\bold{i}\xi S_t/t^{1/\alpha}}]
&=\exp\left(\lambda\int_0^t(\mathrm{e}^{\eta_{c,\alpha,\beta}(\bold{i}\xi F(s)/t^{1/\alpha})}-1)\,\mathrm{d}s\right)
=\exp\left(\lambda\int_0^t(\mathrm{e}^{\frac{1}{t}\eta_{c,\alpha,\beta}(\bold{i}\xi F(s))}-1)\,\mathrm{d}s\right).\label{eq:ftrasfSn}
\end{align}
Again, by the scaling property of the L\'evy exponent $\eta_{c,\alpha,\beta}(\cdot)$ and its definition, we have
\begin{equation}\label{eq:2novsera}
\eta_{c,\alpha,\beta}\left(\bold{i}\xi t^{\frac{1}{2\alpha}}\right)=\eta_{c,\alpha,\beta}\left(\frac{\bold{i}\xi}{(t^{-1/2})^{1/\alpha}}\right)=\sqrt{t}\eta_{c,\alpha,\beta}(\bold{i}\xi),\quad\text{$t>0$, $\xi\in\mathbb{R}$,}
\end{equation}
and 
$\eta_{c\lambda^{1/\alpha},\alpha,\beta}(\bold{i}\xi)=\lambda\eta_{c,\alpha,\beta}(\bold{i}\xi)$  , $\xi\in\mathbb R$.
By these relations and \eqref{eq:ftrasfSn}, for any $t>0$ and $\xi\in\mathbb{R}$, it follows
\begin{align}
\psi_t(\bold{i}\xi)&=\exp\left(\lambda\int_0^t(\mathrm{e}^{\frac{1}{t}\eta_{c,\alpha,\beta}(\bold{i}\xi t^{1/(2\alpha)}F(s))}-1)\,\mathrm{d}s-\lambda\sqrt{t}\eta_{c,\alpha,\beta}(\bold{i}\xi)\right)\nonumber\\
&=\exp\left(\lambda\int_0^t\left(\mathrm{e}^{\frac{1}{\sqrt{t}}\eta_{c,\alpha,\beta}(\bold{i}\xi F(s))}-1-\frac{1}{\sqrt{t}}\eta_{c,\alpha,\beta}(\bold{i}\xi)\right)\,\mathrm{d}s\right)\nonumber\\
&=\exp\left(\lambda\int_0^t\left(\mathrm{e}^{\frac{1}{\sqrt{t}}\eta_{c,\alpha,\beta}(\bold{i}\xi F(s))}-
\mathrm{e}^{\frac{1}{\sqrt{t}}\eta_{c,\alpha,\beta}(\bold{i}\xi)}\right)\,\mathrm{d}s+\lambda t\left(\mathrm{e}^{\frac{1}{\sqrt{t}}\eta_{c,\alpha,\beta}(\bold{i}\xi)}-1-\frac{1}{\sqrt{t}}\eta_{c,\alpha,\beta}(\bold{i}\xi)\right)\right).
\label{eq:8ott1}
\end{align}
By Lemma \ref{le:standardineq} and the definition of the L\'evy exponent $\eta_{c,\alpha,\beta}(\cdot)$, for any $t>0$, $s\in (0,t)$ and 
$\xi\in\mathbb{R}$, we have
\begin{align}
\Big|\mathrm{e}^{\frac{1}{\sqrt{t}}\eta_{c,\alpha,\beta}(\bold{i}\xi F(s))}-
\mathrm{e}^{\frac{1}{\sqrt{t}}\eta_{c,\alpha,\beta}(\bold{i}\xi)}\Big|&\leq\frac{1}{\sqrt t}\Big|\eta_{c,\alpha,\beta}(\bold{i}\xi F(s))-\eta_{c,\alpha,\beta}(\bold{i}\xi)\Big|
\mathrm{e}^{\frac{1}{\sqrt t}\max\{\mathrm{Re}(\eta_{c,\alpha,\beta}(\bold{i}\xi F(s))),\mathrm{Re}(\eta_{c,\alpha,\beta}(\bold{i}\xi))\}}\nonumber\\
&=\frac{1}{\sqrt t}|\eta_{c,\alpha,\beta}(\bold{i}\xi F(s))-\eta_{c,\alpha,\beta}(\bold{i}\xi)|\mathrm{e}^{-\frac{1}{\sqrt t}|c\xi F(s)|^{\alpha}}\nonumber\\
&=\frac{1}{\sqrt t}|1-\bold{i}\beta\tan(\pi\alpha/2)\mathrm{sgn}(\xi)|||c\xi|^{\alpha}-|c\xi F(s)|^{\alpha}|\mathrm{e}^{-\frac{1}{\sqrt t}|c\xi F(s)|^{\alpha}}
\nonumber\\
&=c^{\alpha}\sqrt{1+(\beta\tan(\pi\alpha/2))^2}\frac{(1-F(s)^{\alpha})\mathrm{e}^{-\frac{1}{\sqrt t}|c\xi F(s)|^{\alpha}}}{\sqrt t}|\xi|^{\alpha}.\nonumber
\end{align}
Therefore, for any $t>0$ and $\xi\in\mathbb{R}$,
\begin{equation}\label{eq:17ott3}
\Big|\lambda\int_0^t\left(\mathrm{e}^{\frac{1}{\sqrt{t}}\eta_{c,\alpha,\beta}(\bold{i}\xi F(s))}-
\mathrm{e}^{\frac{1}{\sqrt{t}}\eta_{c,\alpha,\beta}(\bold{i}\xi)}\right)\,\mathrm{d}s\Big|\leq H_1(t)|\xi|^{\alpha}.
\end{equation}
Applying Taylor's formula with integral remainder to the function 
\[
u\mapsto\mathrm{e}^{u\frac{\eta_{c,\alpha,\beta}(\bold{i}\xi)}{\sqrt t}},\quad u\in [0,1],
\]
for any $t>0$ and $\xi\in\mathbb R$, we have
\begin{align}
\mathrm{e}^{\frac{\eta_{c,\alpha,\beta}(\bold{i}\xi)}{\sqrt t}}-1&=\frac{\eta_{c,\alpha,\beta}(\bold{i}\xi)}{\sqrt t}+
\frac{(\eta_{c,\alpha,\beta}(\bold{i}\xi))^2}{t}
\int_0^1(1-u)\mathrm{e}^{u\frac{\eta_{c,\alpha,\beta}(\bold{i}\xi)}{\sqrt t}}\,\mathrm{d}u.\nonumber
\end{align}
Therefore, for any $t>0$ and $\xi\in\mathbb{R}$,
\begin{equation}\label{eq:8ott2}
\Big|\lambda t\left(\mathrm{e}^{\frac{1}{\sqrt{t}}\eta_{c,\alpha,\beta}(\bold{i}\xi)}-1-\frac{1}{\sqrt{t}}\eta_{c,\alpha,\beta}(\bold{i}\xi))
\right)\Big|\leq\frac{\lambda}{2}|\eta_{c,\alpha,\beta}(\bold{i}\xi)|^2=H_2(t)|\xi|^{2\alpha},
\end{equation}
where $H_2(\cdot)$ is the constant function defined in \eqref{eq:K2primo}. 
Note that in \eqref{eq:8ott2} we used the relation
\begin{equation*}
\Big|\mathrm{e}^{\frac{u}{\sqrt t}\eta_{c,\alpha,\beta}(\bold{i}\xi)}\Big|
=\left|\mathrm{e}^{\eta_{c,\alpha,\beta}\left(\frac{\bold{i}\xi}{\left(\frac{\sqrt t}{u}\right)^{1/\alpha}}\right)}\right|
=\mathrm{e}^{-\frac{u}{\sqrt t}|c\xi|^\alpha}\leq 1,\quad u\in [0,1],
\end{equation*}
which follows by \eqref{eq:moduluscarfunc} and Lemma \ref{le:standardineq}.
By Lemma \ref{le:standardineq} and $\max\{\mathrm{Re}z,0\}\leq |z|$, $z\in\mathbb{C}$, we have 
\begin{equation*}
|\mathrm{e}^{z}-1|\leq |z|\mathrm{e}^{|z|},\quad z\in\mathbb C. 
\end{equation*}
The claim follows combining this elementary inequality
with \eqref{eq:8ott1}, \eqref{eq:17ott3} and \eqref{eq:8ott2}.
\\
\noindent{\it Step\,\,2:\,\,Conclusion\,\,of\,\,the\,\,proof.}\\ 
By \eqref{eq:2021Feb15},
for any $t>0$ and $\xi\in\mathbb R$, we have
\begin{align}
\Big|\psi_t\left(\bold{i}\frac{\xi}{t^{1/(2\alpha)}}\right)-1\Big|&\leq\left(\lambda c^{\alpha}\sqrt{1+(\beta\tan(\pi\alpha/2))^2}
\frac{\int_0^t(1-F(s)^{\alpha})\,\mathrm{d}s}{t}|\xi|^{\alpha}+\frac{\lambda c^{2\alpha}}{2t}[1+(\beta\tan(\pi\alpha/2))^2]|\xi|^{2\alpha}\right)\nonumber\\
&\times
\exp\left(\lambda c^{\alpha}\sqrt{1+(\beta\tan(\pi\alpha/2))^2}
\frac{\int_0^t(1-F(s)^{\alpha})\,\mathrm{d}s}{t}|\xi|^{\alpha}+\frac{\lambda c^{2\alpha}}{2t}[1+(\beta\tan(\pi\alpha/2))^2]|\xi|^{2\alpha}\right).\label{eq:9Nov1}
\end{align}
The right-hand side of this latter inequality goes to zero as $t\to\infty$, indeed by de l'Hopital's theorem and the fact that $F(\cdot)$ is a distribution function, we have
\begin{equation*}
\lim_{t\to\infty}\frac{1}{t}\int_0^t(1-F(s)^{\alpha})\,\mathrm{d}s=0.
\end{equation*}
The claim then follows by L\'evy's continuity theorem noticing that, for any $t>0$ and $\xi\in\mathbb R$,
\begin{equation}\label{eq:Nov11mat4}
\psi_{t}\left(\bold{i}\frac{\xi}{t^{1/(2\alpha)}}\right)
=\frac{\mathbb{E}[\mathrm{e}^{\bold{i}(\xi/t^{1/(2\alpha)})t^{-1/(2\alpha)}S_t}]}{\mathrm{e}^{\sqrt{t}\eta_{c\lambda^{1/\alpha},\alpha,\beta}(\bold{i}(\xi/t^{1/(2\alpha)}))}}
=\frac{\mathbb{E}[\mathrm{e}^{\bold{i}\xi(S_t/t^{1/\alpha})}]}{\mathrm{e}^{\eta_{c\lambda^{1/\alpha},\alpha,\beta}(\bold{i}\xi)}},
\end{equation}
where the latter equality follows by \eqref{eq:2novsera}. \\
\noindent$\square$

\noindent{\it Proof\,\,of\,\,Theorem\,\,\ref{prop:speed}.}
By \eqref{eq:2021Feb15}, for any $t>0$ and $\xi\in\mathbb{R}$, we have \eqref{eq:9Nov1}, which, for $\eta\in (0,1]$, we rewrite as
\begin{equation}\label{eq:Nov9II}
\Big|\psi_t\left(\bold{i}\frac{\xi}{t^{1/(2\alpha)}}\right)-1\Big|\leq\left(\frac{H_{1,\eta}(t)}{t^{1-\eta}}|\xi|^{\alpha}
+\frac{H_2}{t}|\xi|^{2\alpha}\right)\exp\left(\frac{H_{1,\eta}(t)}{t^{1-\eta}}|\xi|^{\alpha}
+\frac{H_2}{t}|\xi|^{2\alpha}\right),
\end{equation}
where 
\[
H_{1,\eta}(t):=\lambda c^{\alpha}\sqrt{1+(\beta\tan(\pi\alpha/2))^2}
\frac{\int_0^t(1-F(s)^{\alpha})\,\mathrm{d}s}{t^{\eta}},\quad t>0,
\]
and $H_2(\cdot)\equiv H_2$ is defined by 
\eqref{eq:K2primo}; note that $H_{1,1/2}(\cdot)\equiv H_1(\cdot)$, where $H_1(\cdot)$ is defined by
\eqref{eq:H1primo}.  Letting $\mathbb{P}_X$ denote the law of a random variable $X$, we 
consider the signed measure 
\begin{equation*}
\mu_t(\mathrm{d}x):=\mathbb{P}_{t^{-1/\alpha}S_t}(\mathrm{d}x)-\mathbb{P}_{S}(\mathrm{d}x),\quad t>0,
\end{equation*}
and note that
\[
\widehat{\mu_t}(\xi):=\mathbb{E}[\mathrm{e}^{\bold{i}\xi(S_t/t^{1/\alpha})}]-\mathbb{E}[\mathrm{e}^{\bold{i}\xi S}]
=\mathbb{E}[\mathrm{e}^{\bold{i}\xi S}]\left(\psi_t\left(\bold{i}\frac{\xi}{t^{\frac{1}{2\alpha}}}\right)-1\right),\quad\text{$t>0$, $\xi\in\mathbb R$,}
\]
where the latter equality follows by \eqref{eq:Nov11mat4}. Therefore, by \eqref{eq:Nov9II} and \eqref{eq:moduluscarfunc}, for any $t>0$ and $\xi\in\mathbb{R}\setminus\{0\}$, we have
\begin{align}
\Big|\frac{\widehat{\mu_t}(\xi)}{\xi}\Big|&=\mathrm{e}^{-\lambda|c\xi|^{\alpha}}\frac{1}{|\xi|}\Big|\psi_t\left(\bold{i}\frac{\xi}{t^{\frac{1}{2\alpha}}}\right)-1\Big|\nonumber\\
&\leq\left(\frac{H_{1,\eta}(t)}{t^{1-\eta}}|\xi|^{\alpha-1}+\frac{H_2}{t}|\xi|^{2\alpha-1}\right)
\exp\left(-\lambda|c\xi|^{\alpha}+
\frac{H_{1,\eta}(t)}{t^{1-\eta}}|\xi|^{\alpha}+\frac{H_2}{t}|\xi|^{2\alpha}
\right).\label{eq:mucappello2}
\end{align}
Let $\rho(\cdot)$ be the kernel function provided by Lemma 2.13 in \cite{FMN2}, and for
$\varepsilon>0$ and $x,a\in\mathbb R$, put
$\rho_\varepsilon(x):=\varepsilon^{-1}\rho(x/\varepsilon)$, $f_\varepsilon(x):=\bold{1}_{(-\infty,0]}*\rho_\varepsilon (x)$, and $f_{a,\varepsilon}(x):=f_\varepsilon(x-a)$.
Here the symbol $*$ denotes the convolution product. We continue by first proving the Part $(i)$ of the theorem, and then the Part $(ii)$.
We present a proof which is a little bit more technical than the necessary since we provide a constant involved in the big $O$ notation.\\
\noindent${\it Proof\,\,of\,\,Part\,\,(i).}$\\
We only prove the claim under the assumption \eqref{eq:sigma}. Indeed, if $\eta_0\notin A$, then $\eta_0+\varepsilon\in A$ for any $\varepsilon>0$
sufficiently small.
So let $\eta_0\in A$ and let $t^*>0$ be such that 
\[
\sup_{t\geq t^*}\frac{\int_0^t(1-F(s)^{\alpha})\,\mathrm{d}s}{t^{\eta_0}}\in [0,\infty).
\]
For an arbitrarily fixed
\begin{equation}\label{eq:ineqgamma}
\gamma\in\left(-\frac{\eta_0}{\alpha(1-\eta_0)},1-\frac{\eta_0}{\alpha(1-\eta_0)}\right], 
\end{equation}
we define
\begin{equation}\label{eq:epsilont}
\varepsilon(t):=\frac{1}{K_0 t^{(1-\eta_0)(\gamma+\frac{\eta_0}{\alpha(1-\eta_0)})}},\quad t>0,\quad K_0:=\left(\frac{\lambda}{4H_2}\right)^{\frac{1}{\alpha}}c,
\end{equation}
and consider the family of functions $\{f_{a,\varepsilon(t)}\}_{t>0}$. By Proposition 2.14 in \cite{FMN2}, we have that the functions $f_{a,\varepsilon(t)}$ are smooth test distributions 
in $\mathcal{T}_1(\mathbb R)$ (we refer the reader to Definition 2.9 on p. 11 of \cite{FMN2} for the rigorous definition of such space)
with Fourier transform compactly supported on the interval
\begin{equation}\label{eq:It}
I_t:=[-K_0 t^{(1-\eta_0)(\gamma+\frac{\eta_0}{\alpha(1-\eta_0)})},K_0 t^{(1-\eta_0)(\gamma+\frac{\eta_0}{\alpha(1-\eta_0)})}]
\end{equation}
and $\|\partial f_{a,\varepsilon(t)}\|_{L^1(\mathbb{R},\mathrm{d}x)}=1$.  Note that
\begin{equation}\label{eq:befana}
\alpha\gamma-1\leq 0.
\end{equation}
Indeed, by \eqref{eq:ineqgamma} and the fact that $\alpha\leq 1/(1-\eta_0)$ (since $\eta_0\geq(\alpha-1)/\alpha$), we have
\[
\alpha\gamma-1\leq\alpha-\frac{\eta_0}{1-\eta_0}-1=\alpha-\frac{1}{1-\eta_0}\leq 0.
\]
For $t\geq 1$ and $\xi\in I_t$, we have
\begin{align}
\frac{|\xi|^{2\alpha}}{t}&=|\xi|^{\alpha}\frac{|\xi|^{\alpha}}{t}\leq K_0^\alpha t^{(1-\eta_0)(\alpha\gamma+\frac{\eta_0}{1-\eta_0})-1}|\xi|^{\alpha}\nonumber\\
&=K_0^{\alpha}t^{(1-\eta_0)(\alpha\gamma-1)}|\xi|^{\alpha}
\leq\frac{\lambda}{4H_2}|c\xi|^{\alpha},\label{eq:10ott1}
\end{align}
where for the latter inequality we used that $t^{(1-\eta_0)(\gamma\alpha-1)}\leq 1$ since $t\geq 1$ and 
$\gamma\alpha-1\leq 0$ by \eqref{eq:befana}. Let $t'\geq 1$ be such that
\begin{equation}\label{hp:11ott}
\frac{\int_0^t(1-F(s)^{\alpha})\,\mathrm{d}s}{t}\leq\frac{\lambda c^\alpha}{4\sqrt{2\lambda H_2}}
,\quad\text{for any $t\geq t'$.}
\end{equation} 
Note that such a $t'$ exists since the left-hand side of \eqref{hp:11ott} tends to zero as $t\to\infty$ by de l'Hopital's theorem.
By \eqref{hp:11ott} we have
\begin{equation}\label{eq:11ott4}
\frac{H_{1,\eta_0}(t)}{t^{1-\eta_0}}\leq\lambda c^{\alpha}\sqrt{1+(\beta\tan(\pi\alpha/2))^2}
\frac{\lambda c^\alpha}{4\sqrt{2\lambda H_2}}=\frac{\lambda c^\alpha}{4},\quad\text{for any $t\geq t'$.}
\end{equation}
By \eqref{eq:mucappello2}, \eqref{eq:10ott1} and \eqref{eq:11ott4}, for any $t\geq t'':=\max\{t',t^*\}$ and $\xi\in I_t\setminus\{0\}$,
we have
\begin{equation}\label{eq:mucappello1}
\Big|\frac{\widehat{\mu_t}(\xi)}{\xi}\Big|\leq\left(K'_1\frac{|\xi|^{\alpha-1}}{t^{1-\eta_0}}+H_2\frac{|\xi|^{2\alpha-1}}{t}\right)
\exp\left(-\frac{\lambda|c\xi|^{\alpha}}{2}\right),
\end{equation}
where
\[
K'_1:=\lambda c^{\alpha}\sqrt{1+(\beta\tan(\pi\alpha/2))^2}
\sup_{t\geq t^*}\frac{\int_0^t(1-F(s)^{\alpha})\,\mathrm{d}s}{t^{\eta_0}}.
\]
We note that, for any $\nu,\kappa,\rho>0$,
\begin{align}
\int_{\mathbb R}|\xi|^{\nu-1}\mathrm{e}^{-\kappa|\xi|^{\rho}}\,\mathrm{d}\xi&=2\int_{0}^{\infty}\xi^{\nu-1}\mathrm{e}^{-\kappa\xi^{\rho}}\,\mathrm{d}\xi
=2\int_{0}^{\infty}x^{(\nu-1)/\rho}\mathrm{e}^{-\kappa x}\,\mathrm{d}x^{1/\rho}
=\frac{2\Gamma\left(\frac{\nu}{\rho}\right)}{\rho\kappa^{\nu/\rho}},\label{eq:gamma}
\end{align}
where $\Gamma(\cdot)$ is the Euler gamma function. In particular, for any $t\geq t''$, the function $\xi\mapsto\frac{\widehat{\mu_t}(\xi)}{\xi}$ is integrable. Therefore by Remark 2.11 in \cite{FMN2} we have that the extended Parseval formula applies
and, letting $\widehat{f}(\xi):=\int_{\mathbb R}\mathrm{e}^{\bold{i}\xi x}f(x)\,\mathrm{d}x$ denote the Fourier transform of $f\in L^1(\mathbb R,\mathrm{d}x)$,
we have
\begin{align}
|\mathbb{E}[f_{a,\varepsilon(t)}(S_t/t^{1/\alpha})]-\mathbb{E}[f_{a,\varepsilon(t)}(S)]|
&=\frac{1}{2\pi}\Big|\int_{I_t}\widehat{\partial f_{a,\varepsilon(t)}}(\xi)\frac{\widehat{\mu_t}(-\xi)}{\xi}\,\mathrm{d}\xi\Big|\nonumber\\
&\leq\frac{\|\partial f_{a,\varepsilon(t)}\|_{L^1(\mathbb{R},\mathrm{d}x)}}{2\pi}\int_{I_t}\Big|\frac{\widehat{\mu_t}(-\xi)}{\xi}\Big|\,\mathrm{d}\xi\label{eq:23Feb2021}\\
&=\frac{1}{2\pi}\int_{I_t}\Big|\frac{\widehat{\mu_t}(-\xi)}{\xi}\Big|\,\mathrm{d}\xi,\quad\text{for any $t\geq t''$,}\label{eq:12ott}
\end{align}
where in \eqref{eq:23Feb2021} we used the elementary inequality $|\widehat{f}(\xi)|\leq\|f\|_{L^1(\mathbb{R},\mathrm{d}x)}$, for any $\xi\in\mathbb R$.
By \eqref{eq:mucappello1} and \eqref{eq:gamma}, for any $t\geq t''$, we have
\begin{align}
&\int_{I_t}\Big|\frac{\widehat{\mu_t}(-\xi)}{\xi}\Big|\,\mathrm{d}\xi\leq\frac{K'_1}{t^{1-\eta_0}}
\int_{I_t}|\xi|^{\alpha-1}\mathrm{e}^{-\frac{\lambda|c|^{\alpha}}{2}|\xi|^{\alpha}}\,\mathrm{d}\xi
+\frac{H_2}{t}
\int_{I_t}|\xi|^{2\alpha-1}\mathrm{e}^{-\frac{\lambda|c|^{\alpha}}{2}|\xi|^{\alpha}}\,\mathrm{d}\xi\nonumber\\
&\qquad\qquad\qquad
\leq\frac{K'_1}{t^{1-\eta_0}}\frac{4}{\lambda\alpha c^{\alpha}}+
\frac{H_2}{t}\frac{8}{\alpha\lambda^2 c^{2\alpha}}\leq\frac{4}{\lambda\alpha c^\alpha}\left(K'_1+\frac{2H_2}{\lambda c^\alpha}\right)\frac{1}{t^{1-\eta_0}},\nonumber
\end{align}
where the latter inequality follows noticing that $t''\geq 1$. Combining this with \eqref{eq:12ott}, for any $t\geq t''$, we have
\begin{align}
|\mathbb{E}[f_{a,\varepsilon(t)}(S_t/t^{1/\alpha})]-\mathbb{E}[f_{a,\varepsilon(t)}(S)]|
&\leq\frac{2}{\pi\lambda\alpha c^\alpha}\left(K'_1+\frac{2H_2}{\lambda c^\alpha}\right)\frac{1}{t^{1-\eta_0}}\nonumber\\
&\leq\frac{2}{\pi\lambda\alpha c^\alpha}\left(K'_1+\frac{2H_2}{\lambda c^\alpha}\right)\frac{1}{t^{(1-\eta_0)(\gamma+\frac{\eta_0}{\alpha(1-\eta_0)})}}=B\varepsilon(t),\label{eq:ineq15ott1}
\end{align}
where 
\[
B:=\frac{2K_0}{\pi\lambda\alpha c^\alpha}\left(K'_1+\frac{2H_2}{\lambda c^\alpha}\right)
\]
and the inequality \eqref{eq:ineq15ott1} is a consequence of
$\gamma\leq 1-\frac{\eta_0}{\alpha(1-\eta_0)}$ and $t''\geq 1$. 
By Theorem 2.15 in \cite{FMN2}, noticing that $S$ has a density with respect to the Lebesgue measure bounded above by 
$\frac{1}{\alpha\pi c\lambda^{1/\alpha}}\Gamma\left(\frac{1}{\alpha}\right)$ (see Subsection \ref{subsec:stable}), we have
\[
\mathrm{d}_{\mathrm{Kol}}(S_t/t^{1/\alpha},S)\leq\frac{C}{t^{(1-\eta_0)(\gamma+\frac{\eta_0}{\alpha(1-\eta_0)})}},\quad\text{for any $t\geq t''$,}
\]
where
\begin{equation}\label{eq:C22DIC}
C:=\inf_{\zeta>0}\left\{(1+\zeta)\left(\frac{2}{\pi\lambda\alpha c^\alpha}\left(K'_1+\frac{2H_2}{\lambda c^\alpha}\right)
+\frac{(4H_2)^{\frac{1}{\alpha}}\Gamma\left(\frac{1}{\alpha}\right)}{\pi^{4/3}\alpha c^2\lambda^{2/\alpha}}(4(1+\zeta^{-1})^{1/3}+3^{4/3})\right)\right\}.
\end{equation}
The claim follows noticing that 
\[
\inf_{-\frac{\eta_0}{\alpha(1-\eta_0)}<\gamma\leq1-\frac{\eta_0}{\alpha(1-\eta_0)}}\,
\frac{1}{t^{(1-\eta_0)(\gamma+\frac{\eta_0}{\alpha(1-\eta_0)})}}=\frac{1}{t^{1-\eta_0}},\quad t\geq t''.
\]
\noindent${\it Proof\,\,of\,\,Part\,\,(ii).}$\\
We take
\begin{equation}\label{eq:ineqgamma2}
\gamma\in\left(-\frac{1}{\alpha},\min\left\{2-\frac{1}{\alpha},\frac{1}{\alpha}\right\}\right] 
\end{equation}
and define $\varepsilon(t)$ and $I_t$ as in \eqref{eq:epsilont} and \eqref{eq:It} with $1/2$ in place of $\eta_0$, i.e., 
\begin{equation}\label{eq:epsilontNEW}
\varepsilon(t):=\frac{1}{K_0 t^{\frac{1}{2}(\gamma+\frac{1}{\alpha})}},\quad t>0,
\end{equation}
and $I_t=[-\varepsilon(t)^{-1},\varepsilon(t)^{-1}]$. Then we consider the family of functions $\{f_{a,\varepsilon(t)}\}_{t>0}$ (the function $f_{a,\varepsilon}$ is defined at the beginning of the proof) which, by Proposition 2.14 in \cite{FMN2}, are smooth test distributions 
in $\mathcal{T}_1(\mathbb R)$ with Fourier transform compactly supported on $I_t$. Note that the relations \eqref{eq:10ott1} and \eqref{eq:11ott4} hold with $1/2$ 
in place of $\eta_0$. By the assumption \eqref{eq:sigma2}, it follows
\[
\frac{H_1(t)}{\sqrt t}\leq\frac{K''_1}{t},\quad\text{for any $t>0$,}
\]
where $H_1(\cdot)$ is given by \eqref{eq:H1primo} and
\[
K''_1:=\lambda c^{\alpha}\sqrt{1+(\beta\tan(\pi\alpha/2))^2}
\int_0^\infty(1-F(s)^{\alpha})\,\mathrm{d}s\in [0,\infty).
\]
Combining this with \eqref{eq:mucappello2}, \eqref{eq:10ott1} and \eqref{eq:11ott4} (again with $1/2$ in place of $\eta_0$), 
for any $t\geq t'$ and $\xi\in I_t\setminus\{0\}$, we have
\begin{align}
\Big|\frac{\widehat{\mu_t}(\xi)}{\xi}\Big|
\leq\left(K''_1\frac{|\xi|^{\alpha-1}}{t}+H_2\frac{|\xi|^{2\alpha-1}}{t}\right)
\exp\left(-\frac{\lambda|c\xi|^{\alpha}}{2}\right).\label{eq:15ottpom1}
\end{align}
 As in the proof of Part $(i)$, we have that, for any $t\geq t'$, the function $\xi\mapsto\frac{\widehat{\mu_t}(\xi)}{\xi}$ is integrable, and
one can apply the extended Parseval formula to get \eqref{eq:12ott}, with $\varepsilon(t)$ defined by \eqref{eq:epsilontNEW} and $t'$ in place of $t''$. By \eqref{eq:15ottpom1}
and \eqref{eq:gamma}, for any $t\geq t'$, we have
\begin{align}
&\int_{I_t}\Big|\frac{\widehat{\mu_t}(-\xi)}{\xi}\Big|\,\mathrm{d}\xi\leq\frac{K''_1}{t}
\int_{I_t}|\xi|^{\alpha-1}\mathrm{e}^{-\frac{\lambda|c|^{\alpha}}{2}|\xi|^{\alpha}}\,\mathrm{d}\xi
+\frac{H_2}{t}
\int_{I_t}|\xi|^{2\alpha-1}\mathrm{e}^{-\frac{\lambda|c|^{\alpha}}{2}|\xi|^{\alpha}}\,\mathrm{d}\xi\nonumber\\
&\qquad\qquad\qquad
\leq\frac{K''_1}{t}\frac{4}{\alpha\lambda c^{\alpha}}+
\frac{H_2}{t}\frac{8}{\alpha\lambda^2 c^{2\alpha}}=\frac{4}{\lambda\alpha c^\alpha}\left(K''_1+\frac{2H_2}{\lambda c^\alpha}\right)\frac{1}{t}.\nonumber
\end{align}
Combining this with \eqref{eq:12ott} (with $\varepsilon(t)$ given by \eqref{eq:epsilontNEW} and $t'$ in place of $t''$), for any $t\geq t'$, we have
\begin{align}
|\mathbb{E}[f_{a,\varepsilon(t)}(S_t/t^{1/\alpha})]-\mathbb{E}[f_{a,\varepsilon(t)}(S)]|
&\leq\frac{2}{\pi\lambda\alpha c^\alpha}\left(K''_1+\frac{2H_2}{\lambda c^\alpha}\right)\frac{1}{t}\nonumber\\
&\leq\frac{2}{\pi\lambda\alpha c^\alpha}\left(K''_1+\frac{2H_2}{\lambda c^\alpha}\right)\frac{1}{t^{\frac{1}{2}(\gamma+\frac{1}{\alpha})}}=B'\varepsilon(t),\label{eq:15ottpom2}
\end{align}
where 
\[
B':=\frac{2K_0}{\pi\lambda\alpha c^\alpha}\left(K''_1+\frac{2H_2}{\lambda c^\alpha}\right)
\]
and the inequality \eqref{eq:15ottpom2} is a consequence of
$\gamma\leq 2-\frac{1}{\alpha}$ and $t'\geq 1$. By Theorem 2.15 in \cite{FMN2}, noticing that $S$ has a density (with respect to the Lebesgue measure) bounded above by 
$\frac{1}{\alpha\pi c\lambda^{1/\alpha}}\Gamma\left(\frac{1}{\alpha}\right)$, we have
\[
\mathrm{d}_{\mathrm{Kol}}(S_t/t^{1/\alpha},S)\leq\frac{C'}{t^{\frac{1}{2}(\gamma+\frac{1}{\alpha})}},\quad\text{for any $t\geq t'$,}
\]
where the constant $C'$ is defined as $C$ in \eqref{eq:C22DIC}, but with $K''_1$ in place of $K'_1$.
The claim follows taking the infimum over $\gamma$, which satisfies the constraint 
\eqref{eq:ineqgamma2}. Indeed, if $\alpha\in (1,2]$, then $1/\alpha<2-1/\alpha$, and so
\[
\inf_{-1/\alpha<\gamma\leq\min\left\{1/\alpha,2-1/\alpha\right\}}\frac{1}{t^{\frac{1}{2}(\gamma+\frac{1}{\alpha})}}=t^{-\frac{1}{\alpha}},\quad t\geq t';
\]
if, instead, $\alpha\in (0,1]$, then $1/\alpha\geq 2-1/\alpha$, and so
\[
\inf_{-1/\alpha<\gamma\leq\min\left\{1/\alpha,2-1/\alpha\right\}}\frac{1}{t^{\frac{1}{2}(\gamma+\frac{1}{\alpha})}}=t^{-1},\quad t\geq t'.
\]
\noindent$\square$

\subsection{Proof of Lemma \ref{le:standardineq}}

For ease of notation, let $a_i:=\mathrm{Re}z_i$ and $b_i:=\mathrm{Im}z_i$, $i=1,2$.
Hereon, without loss of generality, we assume $a_1\geq a_2$. We have:
\begin{align}
|\re^{z_1}-\re^{z_2}|
=|\re^{a_1+\im b_1}-\re^{a_2+\im b_1}+\re^{a_2+\im b_1}-\re^{a_2+\im b_2}|
=|\re^{a_1}-\re^{a_2}+\re^{a_2}(1-\re^{\im(b_2-b_1)})|.\nonumber
\end{align}
By the mean value theorem we have $\re^{a_1}-\re^{a_2}=(a_1-a_2)\re^c$, for some $c=c(a_1,a_2)\in [a_2,a_1]$. Therefore
\[
|\re^{z_1}-\re^{z_2}|=|(a_1-a_2)\re^{c}+\re^{a_2}(1 -\re^{\im(b_2-b_1)}) |\leq |(a_1-a_2)\re^{a_1}+\re^{a_2}(1-\re^{\im(b_2-b_1)})|.
\]
Setting $\alpha:=a_1-a_2\geq 0$ and $\beta:= b_2-b_1$, we have
\[
|(a_1-a_2)\re^{a_1}+\re^{a_2}(1-\re^{\im(b_2-b_1)})|=\re^{a_1}|\alpha+\re^{-\alpha}(1-\re^{\im\beta})|=
\re^{a_1}|\alpha+\re^{-\alpha}(1-\cos\beta-\im\sin\beta)|.
\]
The claim follows if we prove the inequality
\begin{equation}\label{eq:finalineq}
|\alpha+\re^{-\alpha}(1-\cos\beta)-\im\re^{-\alpha}\sin\beta|\leq\sqrt{\alpha^2+\beta^2}.
\end{equation}
Indeed $|z_1-z_2|=\sqrt{\alpha^2+\beta^2}$ and
\[
|\alpha+\re^{-\alpha}(1-\cos\beta-\im\sin\beta)|=|\alpha+\re^{-\alpha}(1-\cos\beta)-\im\re^{-\alpha}\sin\beta|.
\]
Note that \eqref{eq:finalineq} is equivalent to
\begin{align}
f_1(\alpha,\beta):=[\alpha+\re^{-\alpha}(1-\cos\beta)]^{2}+\re^{-2\alpha}(\sin\beta)^2
\leq\alpha^2+\beta^2=:f_2(\alpha,\beta).\label{eq:fine1}
\end{align}
To prove this latter inequality,  we start noticing that $f_1(\alpha,0)=f_2(\alpha,0)=\alpha^2$ and
\begin{align}
f_1(\alpha,\beta)&=\alpha^2+2\alpha\re^{-\alpha}(1-\cos\beta)+2\re^{-2\alpha}(1-\cos\beta)
=\alpha^2+2(\alpha+\mathrm{e}^{-\alpha})\mathrm{e}^{-\alpha}(1-\cos\beta).\nonumber
\end{align}
Then we distinguish two cases: $\beta>0$ and
$\beta<0$.  Assume first $\beta>0$. For any $\gamma>0$, we have
\[
\frac{\partial f_1(\alpha, \gamma)}{\partial\gamma}=2(\alpha+\re^{-\alpha})\re^{-\alpha}\sin\gamma\leq
2\gamma=\frac{\partial f_2(\alpha,\gamma)}{\partial\gamma},
\]
where we used the elementary inequalities $\sin\gamma\leq|\sin\gamma|\leq|\gamma|=\gamma$ and $0<(\alpha+\re^{-\alpha})\re^{-\alpha}\leq 1$ (the latter follows
by the elementary relations $\re^{\alpha}\geq1+\alpha\geq\mathrm{e}^{-\alpha}+\alpha$; recall that $\alpha\geq 0$). Therefore
\begin{equation}\label{eq:fine}
f_1(\alpha,\beta)=f_1(\alpha, 0)+\int_{0}^{\beta}\frac{\partial f_1(\alpha,\gamma)}{\partial\gamma}\,\mathrm{d}\gamma\leq
f_2(\alpha,0)+\int_{0}^{\beta}\frac{\partial f_2(\alpha,\gamma)}{\partial\gamma}\,\mathrm{d}\gamma=f_2(\alpha,\beta),
\end{equation}
and the claim \eqref{eq:fine1} for $\beta>0$ is proved. Now, assume $\beta<0$. By \eqref{eq:fine} we have
$f_1(\alpha,-\beta)\leq f_2(\alpha,-\beta)$ and the claim \eqref{eq:fine1} for $\beta<0$ follows noticing that $f_i(\alpha,-\beta)=f_i(\alpha,\beta)$, $i=1,2$.
\\
\\
\noindent$\bold{Acknowledgments}$
We would like to thank the Editor and two anonymous Reviewers for a careful reading of the paper.
We would like to acknowledge support for the project titled “Epidemics and Counting Structures in Erd\"os–R\'enyi Random Graphs” from the Istituto Nazionale di Alta Matematica “Francesco Severi”.

\end{document}